\documentclass[10pt,twosided]{article}
\usepackage[tbtags]{amsmath}
\usepackage{epsfig,amstext,amssymb,amsthm,latexsym}
\allowdisplaybreaks[4]
\usepackage{epsfig}
\usepackage{color}
\usepackage{graphicx, subfigure}

\usepackage{makeidx}
\usepackage{multirow}
\usepackage{multicol}
\usepackage[dvipsnames,svgnames,table]{xcolor}
\usepackage{graphicx}
\usepackage{epstopdf}
\usepackage{ulem}
\usepackage{hyperref}
\usepackage{amsmath}

\usepackage{amssymb,color}

\definecolor{c20}{rgb}{0.,0.7,0.}
\definecolor{c30}{rgb}{0.,0.,1.}
\definecolor{c40}{rgb}{1,0.1,0.7}
\definecolor{c50}{rgb}{1,0,0}
\definecolor{c60}{rgb}{1,0.9,0.1}

\newtheorem{theorem}{Theorem}[section]

\newtheorem{lemma}{Lemma}[section]

\newtheorem{remark}{Remark}[section]

\numberwithin{equation}{section}

\allowdisplaybreaks[4]

\numberwithin{equation}{section}

\newcommand{\prooftheo}[1]{ \textsc{Proof of Theorem} \ref{#1} }

\newcommand{\prooflem}[1]{\textsc{Proof of Lemma} \ref{#1}}

\newcommand{\pk}[1]{\mathbb{P} \left( #1 \right) }

\newcommand{\COM}[1]{}

\newcommand{\R}{\mathbb{R}}
\newcommand{\inr}{\in \R}

\begin{document}

\title{Second order expansions of distributions of maxima of bivariate Gaussian triangular arrays under power
normalization
\thanks{The work was supported by the National Natural Science Foundation of China (grant No. 11501113 and No. 11601330) and the Key Project of Fujian
Education Committee (grant No. JA15045)}}
\author{Zhichao Weng$^a$ \qquad Xin Liao$^b$\thanks{Corresponding author. Email: liaoxin2010@163.com}
\\{\small $^a$School of Economics and Management, Fuzhou University, Fujian, 350116, China}\\
{\small $^b$Business School, University of Shanghai for Science and Technology, Shanghai, 200093, China}}

\maketitle

{\bf Abstract.} In this paper, we study second order expansions of distributions of maxima of bivariate Gaussian triangular arrays under power normalization. Numerical
analysis are given to compare the asymptotic behaviors under power normalization with the asymptotic behaviors under linear normalization derived by Hashorva et al. (2016).

{\bf keywords:} Second order expansion; Maximum; Bivariate Gaussian triangular array; Power normalization

\section{Introduction}\label{sec1}
 Let $\{(X_{n,k},Y_{n,k}), 1\le k\le n, n\ge 1\}$ be a
triangular array of independent standard bivariate Gaussian random
vector with correlations $\{\rho_n,n\ge 1\}$ and joint distribution
function $F$. H\"{u}sler and Reiss (1989) considered the asymptotic
behavior of distribution of maxima with correlation coefficient
varying as the sample size increases. Under the so-called
H\"{u}sler-Reiss condition
\begin{equation}\label{eq1.3}
\lim_{n\to \infty}\frac{1}{2}b_n^2(1-\rho_n)=\lambda^2
\end{equation}
with $\lambda\in [0,\infty]$, H\"{u}sler and Reiss (1989) showed that
$$\lim_{n\to \infty}\sup_{x,y\inr}|F^n(b_n+x/b_n,b_n+y/b_n)-H_{\lambda}(x,y)|=0$$
holds, where the norming constant $b_n$ is given by
\begin{equation}\label{eq1.2}
n(1-\Phi(b_n))=1,
\end{equation}
with $\Phi$ standing for the standard Gaussian distribution, and the
max-stable H\"{u}sler-Reiss distribution $H_{\lambda}$ is given by
$$H_{\lambda}(x,y)=\exp\left(-\Phi(\lambda+\frac{x-y}{2\lambda})e^{-y}-\Phi(\lambda+\frac{y-x}{2\lambda})e^{-x}\right), \quad x,y\inr$$
with $H_0(x,y)=\exp\left(-e^{-\min(x,y)}\right)$ and
$H_{\infty}(x,y)=\Lambda(x)\Lambda(y)$  with
$\Lambda(x)=\exp\left(-e^{-x}\right)$, $x\inr$.

Extensions the work of H\"{u}sler and Reiss (1989) can be found in
recent literature. Hashorva (2005, 2006) considered the case of
triangular arrays of independent elliptical random vectors, Hashorva
and Ling (2016) extended the results to the case of bivariate
skew-normal triangular array. Motivated by the work of Nair (1981)
and Frick and Reiss (2013), for the H\"{u}ser-Reiss model Hashorva
et al. (2016) established the higher-order expansions of
distributions of maxima under the refined H\"{u}sler-Reiss
conditions and Liao and Peng (2014) considered its associated
uniform convergence rates.

In this paper, we are interested in the rate of convergence of the
distribution of maxima of H\"{u}sler-Reiss model under power
normalization. For univariate case, it's well known that
$$\lim_{n\to \infty}\sup_{x\inr}|\Phi^n(b_n+x/b_n)-\Lambda(x)|=0$$
holds with $b_n$ given by \eqref{eq1.2}, cf. Resnick (1987) and Nair
(1981). In view of Mohan and Ravi (1993), let $\alpha_n=b_n$ and
$\beta_n=b_n^{-2}$, we have
\begin{equation}\label{eq1.4}
\lim_{n\to \infty}\Phi^n(b_nx^{b_n^{-2}})=\Phi_{1}(x),
\end{equation}
where $\Phi_{1}(x)=\exp\left(-x^{-1}\right), x>0$, one of six-type power-stable
distributions given by Pancheva(1985). For recent work on maxima
under power normalization, see Mohan and Subramanya (1991), Mohan and
Ravi (1993), Subramanya (1994), Barakat et al. (2010) and Peng et
al. (2013). In this paper we will show that under the
H\"{u}sler-Reiss condition \eqref{eq1.3}
\begin{equation}\label{eq1.5}
\lim_{n\to \infty}\sup_{x,y>0}|F^n(b_nx^{b_n^{-2}},b_ny^{b_n^{-2}})-H_{\lambda}(\ln x, \ln y)|=0
\end{equation}
holds with $H_0(\ln x,\ln y)=\exp\left(-(\min(x,y))^{-1}\right)$ and
$H_{\infty}(\ln x,\ln y)=\exp\left(-x^{-1}\right)\exp\left(-y^{-1}\right)$. Furthermore,  the
rate of convergence in \eqref{eq1.4} and \eqref{eq1.5} will be
investigated, respectively.

The rest of the paper is organized as follows. In Section \ref{sec2}
we present the main results. Numerate analysis provided in Section \ref{secN}
compare the asymptotic behaviors under power and linear normalization.
All the proof are relegated to Section
\ref{sec3}.

\section{Main results}\label{sec2}

In this section, we provide the main results with respect to
limiting distribution of maxima under power normalization under the
H\"{u}sler-Reiss condition \eqref{eq1.3} and its second-order
expansions providing some refined H\"{u}sler-Reiss condition hold.
In the following we shall denote throughout by $b_n$ the constants
defined in \eqref{eq1.2}. First we state \eqref{eq1.5} as the
following result.
\begin{theorem}\label{thm2.0}
For the considered  bivariate Gaussian triangular array, the
H\"{u}sler-Reiss condition \eqref{eq1.3} holds if and only if
\eqref{eq1.5} holds.
\end{theorem}

To establish the higher-order expansion of the distribution of
maxima in H\"{u}sler-Reiss model, we need to refine the
H\"{u}sler-Reiss condition \eqref{eq1.3}. There are three cases to
be considered, i.e., $\lambda\in (0,\infty)$, $\lambda=0$ and
$\lambda=\infty$, respectively.

For $\lambda\in (0,\infty)$, the result is given as follows.

\begin{theorem}\label{th2.1}
If the second-order H\"{u}sler-Reiss condition
\begin{equation}\label{eq2.1}
\lim_{n\to \infty}b_n^2(\lambda-\lambda_n)=\tau\inr
\end{equation}
holds with $\lambda_n=(\frac{1}{2}b_n^2(1-\rho_n))^{1/2}$ and $\lambda\in (0,\infty)$, then
for $x>0$, $y>0$ we have
\begin{eqnarray}\label{addeq2.2}
\lim_{n\to \infty}b_n^2\left(F^n(b_nx^{b_n^{-2}},b_ny^{b_n^{-2}})-H_{\lambda}(\ln x, \ln y)\right)
=\kappa(x,y,\lambda,\tau)H_{\lambda}(\ln x,\ln y)
\end{eqnarray}
with $\kappa(x,y,\lambda,\tau)$ given by \eqref{addeq3.12}.
\end{theorem}

\begin{remark}\label{rem2.2}
(i) Let $\gamma_n=(\lambda-\lambda_n)^{-1}$. If \eqref{eq2.1} does not converge but
$\gamma_n$ and $b_n^2$ are the same order,
then $b_{n}^{-2}$
and $ F^n\left( b_nx^{b_n^{-2}},b_ny^{b_n^{-2}} \right) -H_{\lambda}(\ln x,\ln y)$ are the same order.

(ii) If $\lim_{n\to \infty} b_n^2/\gamma_n=\pm\infty$, with arguments similar to
that of Theorem 2.2, we can show that
\begin{eqnarray}\label{addeq2.3}
\lim_{n\to \infty} \gamma_n \left( F^n\left( b_nx^{b_n^{-2}},b_ny^{b_n^{-2}} \right) -H_{\lambda}(\ln x,\ln y)\right)
= 2x^{-1}\varphi\left( \lambda+ \frac{\ln\frac{y}{x}}{2\lambda} \right)H_{\lambda}(\ln x,\ln y).
\end{eqnarray}

(iii) Conversely, for the bivariate Gaussian triangular arrays with correlations $\{\rho_n\}$ satisfying \eqref{eq1.3}, we have the following assertions under power normalization:
(a) if \eqref{addeq2.2} holds, then \eqref{eq2.1} holds.
(b) if $ F^n\left( b_nx^{b_n^{-2}},b_ny^{b_n^{-2}} \right) -H_{\lambda}(\ln x,\ln y)$ and $b_n^{-2}$ are the same order,
then $\gamma_n$ and $b_{n}^2$ are the same order.
(c) if \eqref{addeq2.3} holds, then $\lim_{n\to \infty} b_n^2/\gamma_n=\pm \infty$.
(d) if $ F^n\left( b_nx^{b_n^{-2}},b_ny^{b_n^{-2}} \right) -H_{\lambda}(\ln x,\ln y)$ and $\gamma_n$ are the same order, then $\gamma_n$ and $b_n^2$ are the same order.
\end{remark}

\begin{remark}\label{rem2.3}
For the case of $\lambda \in (0,\infty)$, if \eqref{eq2.1} does not converge, and $\gamma_n$ and $b_n^2$ are not
the same order, there may be no convergence rates for the extremes. An example is: suppose that the bivariate
Gaussian triangular arrays have correlations $\{\rho_n \}$ satisfying \eqref{eq1.3}. Furthermore, assume that
$\lim_{n\to \infty} b_{2n}^2/\gamma_{2n}=0$ and
$\lim_{n\to \infty} b_{2n+1}^2/\gamma_{2n+1}=\infty$. Hence, by Theorem \ref{th2.1} and Remark \ref{rem2.2} (ii),
we have
\begin{eqnarray*}
&&\lim_{n\to \infty} b_{2n}^2\left( F^{2n}\left( b_{2n}x^{b_{2n}^{-2}},b_{2n}y^{b_{2n}^{-2}} \right) -H_{\lambda}(\ln x,\ln y)\right)
=\kappa(x,y,\lambda,0)H_{\lambda}(\ln x,\ln y)
\end{eqnarray*}
and
\begin{eqnarray*}
\lim_{n\to \infty} \gamma_{2n+1}\left( F^{2n+1}\left( b_{2n+1}x^{b_{2n+1}^{-2}},b_{2n+1}y^{b_{2n+1}^{-2}} \right) -H_{\lambda}(\ln x,\ln y)\right)
= 2x^{-1} \varphi\left( \lambda+\frac{\ln \frac{y}{x}}{2\lambda} \right) H_{\lambda}(\ln x,\ln y).
\end{eqnarray*}
\end{remark}

Next we give the results of two extreme cases: $\lambda=\infty$ and $\lambda=0$ with
different refined conditions. The following theorem considers the case of $\lambda=\infty$.

\begin{theorem}\label{th2.2}
For $\rho_n\in [-1,1)$, assume that $\lim_{n\to \infty} \frac{\ln b_n}{b_n^2(1-\rho_n)}=0$.
Then for all $x>0$, $y>0$ we have
\begin{equation}\label{eq2.2}
\lim_{n\to \infty}b_n^2\left(F^n(b_nx^{b_n^{-2}},b_ny^{b_n^{-2}})-H_{\infty}(\ln x, \ln y)\right)
=(s(x)+s(y))H_{\infty}(\ln x,\ln y)
\end{equation}
with $s(x)$ given by \eqref{eq3.1}.
\end{theorem}
\COM{
\begin{remark}
For the univariate case, we have (\pzx{only guess here!!})
\[
\lim_{n\to
\infty}b_n^2\left(\Phi^{n}(b_nx^{b_n^{-2}})-\Phi_{1}(x)\right)
=s(x))\Phi_{1}(x)
\]
with $s(x)$ given by \eqref{eq3.1}
\end{remark}
}

For the case of $\lambda=0$, we have the following result.

\begin{theorem}\label{th2.3}
For $\rho_n \in (0,1]$, assume that $\lim_{n\to \infty}b_{n}^6(1-\rho_n)=0$. Then for
$x>0$, $y>0$ we have
\begin{equation}\label{eq2.3}
\lim_{n\to \infty}b_n^2\left(F^n(b_nx^{b_n^{-2}},b_ny^{b_n^{-2}})-H_{0}(\ln x, \ln y)\right)
=s(\min(x,y))H_{0}(\ln x,\ln y)
\end{equation}
with $s(x)$ given by \eqref{eq3.1}.
\end{theorem}

\section{Numerical analysis}\label{secN}
In this section, numerical studies are presented to illustrate the accuracy of second order expansions of
$F^n$ under two different normalization, i.e., the finite behaviors under power normalization derived in this
paper and that under linear normalization given by Hashorva et al. (2016). We shall discuss three particular cases:\\
$(i)$ $\lambda\in (0,\infty)$ with
\begin{equation}\label{addeq1}
\rho_n=1-\frac{2\lambda^2}{b_n^2}+\frac{4\tau
\lambda}{b_n^4}-\frac{2\tau^2}{b_n^6},
\end{equation}
where $b_n$ satisfies \eqref{eq1.2}, which implies that  condition \eqref{eq2.1} holds;\\
$(ii)$ $\rho_{n}=-1,0$ implying $\lambda=\infty$;\\
 $(iii)$ $\rho_{n}=1$ implying
$\lambda=0$.

For the case of power normalization, we calculate the
actual values $F^n(b_nx^{b_n^{-2}},b_ny^{b_n^{-2}})$, the first-order
asymptotics $L_1^{p}=H_{\lambda}(\ln x,\ln y)$, the second-order asymptotics according to the values of $\rho_{n}$ with finite $n$, i.e.,
\begin{itemize}
\item[$(i)$.] if $\rho_{n}$ is given by \eqref{addeq1} with
fixed $\lambda$ and $\tau$, then in view of \eqref{addeq2.2}  the
second-order asymptotics are given by
$L_2^{p}=H_{\lambda}(\ln x,\ln y)\left(1+b_n^{-2}\kappa(x,y,\lambda,\tau)\right)$;
\item[ $(ii)$.] if $\rho_{n}=-1,0$, by \eqref{eq2.2} the second-order
 asymptotics are given by
$L_3^{p}=H_{\infty}(\ln x,\ln y)\left(1+b_{n}^{-2}(s(x)+s(y))\right)$; and
\item[$(iii)$.] if $\rho_{n}=1,n\ge 1$ by \eqref{eq2.3}  the second-order
asymptotics are given by
$L_4^{p}=H_{0}(\ln x,\ln y)\left(1+b_{n}^{-2}s(\min(x,y))\right)$.
\end{itemize}

For the linear normalization case, by Theorem 2.1 and Theorem 2.3 in Hashorva et al. (2016), we have
the first-order
asymptotics $L_1^{l}=H_{\lambda}(x, y)$, the second-order asymptotics according to the values of $\rho_{n}$ with finite $n$, i.e.,
\begin{itemize}
\item[$(i)$.] if $\rho_{n}$ is given by \eqref{addeq1} with
fixed $\lambda$ and $\tau$, then the
second-order asymptotics are given by
$L_2^{l}=H_{\lambda}(x,y)(1+b_n^{-2}\iota(x,y,\lambda,\tau))$;
\item[ $(ii)$.] if $\rho_{n}=-1,0$, the second-order
 asymptotics are given by
$L_3^{l}=H_{\infty}(x,y)\left(1+b_{n}^{-2}(\omega(x)+\omega(y))\right)$; and
\item[$(iii)$.] if $\rho_{n}=1,n\ge 1$ the second-order
asymptotics are given by
$L_4^{l}=H_{0}(x,y)\left(1+b_{n}^{-2}\omega(\min(x,y))\right)$,
\end{itemize}
where $\omega(x)$ and $\iota(x,y,\lambda,\tau)$ are given by
$$\omega(x)=2^{-1}(x^2+2x)e^{-x},$$
$$\iota(x,y,\lambda,\tau)=\omega(x)\Phi(\lambda+\frac{y-x}{2\lambda})
+\omega(y)\Phi(\lambda+\frac{x-y}{2\lambda})+
(2\tau-\lambda(\lambda^2+x+y+2))e^{-x}\varphi(\lambda+\frac{y-x}{2\lambda}).$$

To compare the accuracy of actual values with its asymptotics, let $\Delta_i^{p}=|F^n(b_nx^{b_n^{-2}},b_ny^{b_n^{-2}})-L_i^{p}|$ and $\Delta_i^{l}=|F^n(b_n+x/b_n,b_n+y/b_n)-L_i^{l}|$, $i=1,2,3,4$ denote the absolute
errors under power and linear normalization, respectively.
We use $\mathbf{R}$ to calculate the absolute errors with sample sizes $n=1000$ and $n=10000$,
and fixed $\lambda$, $\tau$, which are documented Table \ref{table:1}-\ref{table:4}.
These tables show that accuracies of the first and the second order asymptotics under two different normalization can be improved as $n$ becomes large.

In order to show the accuracy of all asymptotics with varying $x$ and $y$, we plot the
actual values and its asymptotics with fixed $\lambda$, $\tau$ and $n=10^3$ by using $\mathbf{R}$. Power normalization cases are illustrated in Figures 1 and 2, where Figure 1 compares the actual values with above three asymptotics with $x=y$, Figure 2 compares the difference of the actual value with above three asymptotics by contour line in the plane. The cases of linear normalization are illustrated in Figures 1 and 2 in Hashorva et al. (2016).

According to Figures 1-2 in Hashorva et al. (2016), Figures 1-2 and Tables \ref{table:1}-\ref{table:4}, we have the following findings: i) The asymptotics
under linear normalization are more closer to its actual values except small $x$.
ii) Under two different normalization, the second
order asymptotics are closer to the actual values as small $x$ except few special cases, contrary to the case of large
$x$, which shows that the first order asymptotics may be better.

\begin{table}\centering
\caption{\c{Absolute errors between actual values and their asymptotics for the case of $\lambda=2$, $\tau=3$}}
\label{table:1}

\vspace{0.5cm}
\begin{tabular}{|p{41pt}|p{32pt}|p{38pt}|p{32pt}|p{38pt}|p{32pt}|p{38pt}|p{32pt}|p{38pt}|}
\hline
\parbox{41pt}{\raggedright \multirow{2}{*}{
(x,y)
}} & \multicolumn{4}{|l|}{\parbox{140pt}{\raggedright
n=1000
}} & \multicolumn{4}{|l|}{\parbox{140pt}{\raggedright
n=10000
}} \\
\cline{2-9}
 & \parbox{32pt}{\raggedright
$\Delta_1^{p}$
} & \parbox{38pt}{\raggedright
$\Delta_1^{l}$
} & \parbox{32pt}{\raggedright
$\Delta_2^{p}$
} & \parbox{38pt}{\raggedright
$\Delta_2^{l}$
} & \parbox{32pt}{\raggedright
$\Delta_1^{p}$
} & \parbox{38pt}{\raggedright
$\Delta_1^{l}$
} & \parbox{32pt}{\raggedright
$\Delta_2^{p}$
} & \parbox{38pt}{\raggedright
$\Delta_2^{l}$
} \\
\hline
\parbox{41pt}{\raggedright
(0.5,0.5)
} & \parbox{32pt}{\raggedright
0.00133
} & \parbox{38pt}{\raggedright
0.01646
} & \parbox{32pt}{\raggedright
0.00078
} & \parbox{38pt}{\raggedright
0.00114
} & \parbox{32pt}{\raggedright
0.00106
} & \parbox{38pt}{\raggedright
0.01128
} & \parbox{32pt}{\raggedright
0.00039
} & \parbox{38pt}{\raggedright
0.0007
} \\
\hline
\parbox{41pt}{\raggedright
(1,1)
} & \parbox{32pt}{\raggedright
0.00239
} & \parbox{38pt}{\raggedright
0.04272
} & \parbox{32pt}{\raggedright
0.00241
} & \parbox{38pt}{\raggedright
0.00217
} & \parbox{32pt}{\raggedright
0.00197
} & \parbox{38pt}{\raggedright
0.03001
} & \parbox{32pt}{\raggedright
0.00135
} & \parbox{38pt}{\raggedright
0.00098
} \\
\hline
\parbox{41pt}{\raggedright
(1,0.5)
} & \parbox{32pt}{\raggedright
0.00233
} & \parbox{38pt}{\raggedright
0.02723
} & \parbox{32pt}{\raggedright
0.00192
} & \parbox{38pt}{\raggedright
0.00022
} & \parbox{32pt}{\raggedright
0.00186
} & \parbox{38pt}{\raggedright
0.01897
} & \parbox{32pt}{\raggedright
0.00108
} & \parbox{38pt}{\raggedright
0.00001
} \\
\hline
\parbox{41pt}{\raggedright
(2,1)
} & \parbox{32pt}{\raggedright
0.00901
} & \parbox{38pt}{\raggedright
0.05237
} & \parbox{32pt}{\raggedright
0.00874
} & \parbox{38pt}{\raggedright
0.00675
} & \parbox{32pt}{\raggedright
0.00598
} & \parbox{38pt}{\raggedright
0.03752
} & \parbox{32pt}{\raggedright
0.00579
} & \parbox{38pt}{\raggedright
0.0033
} \\
\hline
\parbox{41pt}{\raggedright
(3,3)
} & \parbox{32pt}{\raggedright
0.07201
} & \parbox{38pt}{\raggedright
0.04517
} & \parbox{32pt}{\raggedright
0.02066
} & \parbox{38pt}{\raggedright
0.01957
} & \parbox{32pt}{\raggedright
0.04978
} & \parbox{38pt}{\raggedright
0.03459
} & \parbox{32pt}{\raggedright
0.01432
} & \parbox{38pt}{\raggedright
0.01011
} \\
\hline
\parbox{41pt}{\raggedright
(3,5)
} & \parbox{32pt}{\raggedright
0.08486
} & \parbox{38pt}{\raggedright
0.0289
} & \parbox{32pt}{\raggedright
0.01672
} & \parbox{38pt}{\raggedright
0.01684
} & \parbox{32pt}{\raggedright
0.05934
} & \parbox{38pt}{\raggedright
0.02264
} & \parbox{32pt}{\raggedright
0.01229
} & \parbox{38pt}{\raggedright
0.00895
} \\
\hline
\parbox{41pt}{\raggedright
(2,3)
} & \parbox{32pt}{\raggedright
0.05262
} & \parbox{38pt}{\raggedright
0.05552
} & \parbox{32pt}{\raggedright
0.01079
} & \parbox{38pt}{\raggedright
0.01649
} & \parbox{32pt}{\raggedright
0.03623
} & \parbox{38pt}{\raggedright
0.04132
} & \parbox{32pt}{\raggedright
0.00735
} & \parbox{38pt}{\raggedright
0.00839
} \\
\hline
\parbox{41pt}{\raggedright
(5,5)
} & \parbox{32pt}{\raggedright
0.09987
} & \parbox{38pt}{\raggedright
0.01059
} & \parbox{32pt}{\raggedright
0.03115
} & \parbox{38pt}{\raggedright
0.01225
} & \parbox{32pt}{\raggedright
0.07058
} & \parbox{38pt}{\raggedright
0.00899
} & \parbox{32pt}{\raggedright
0.02314
} & \parbox{38pt}{\raggedright
0.00679
} \\
\hline
\parbox{41pt}{\raggedright
(5,9)
} & \parbox{32pt}{\raggedright
0.10039
} & \parbox{38pt}{\raggedright
0.00556
} & \parbox{32pt}{\raggedright
0.02004
} & \parbox{38pt}{\raggedright
0.00712
} & \parbox{32pt}{\raggedright
0.07231
} & \parbox{38pt}{\raggedright
0.00475
} & \parbox{32pt}{\raggedright
0.01684
} & \parbox{38pt}{\raggedright
0.004
} \\
\hline
\parbox{41pt}{\raggedright
(10,10)
} & \parbox{32pt}{\raggedright
0.09729
} & \parbox{38pt}{\raggedright
0.00009
} & \parbox{32pt}{\raggedright
0.04002
} & \parbox{38pt}{\raggedright
0.00046
} & \parbox{32pt}{\raggedright
0.07223
} & \parbox{38pt}{\raggedright
0.00009
} & \parbox{32pt}{\raggedright
0.03268
} & \parbox{38pt}{\raggedright
0.00029
} \\
\hline
\parbox{41pt}{\raggedright
(10,20)
} & \parbox{32pt}{\raggedright
0.08437
} & \parbox{38pt}{\raggedright
0.00005
} & \parbox{32pt}{\raggedright
0.02078
} & \parbox{38pt}{\raggedright
0.00024
} & \parbox{32pt}{\raggedright
0.06431
} & \parbox{38pt}{\raggedright
0.00004
} & \parbox{32pt}{\raggedright
0.02041
} & \parbox{38pt}{\raggedright
0.00015
} \\
\hline
\parbox{41pt}{\raggedright
(7,10)
} & \parbox{32pt}{\raggedright
0.10072
} & \parbox{38pt}{\raggedright
0.00091
} & \parbox{32pt}{\raggedright
0.02716
} & \parbox{38pt}{\raggedright
0.00231
} & \parbox{32pt}{\raggedright
0.07369
} & \parbox{38pt}{\raggedright
0.00084
} & \parbox{32pt}{\raggedright
0.02289
} & \parbox{38pt}{\raggedright
0.00138
} \\
\hline
\parbox{41pt}{\raggedright
(20,20)
} & \parbox{32pt}{\raggedright
0.06994
} & \parbox{38pt}{\raggedright
\small{4.02$\times10^{-9}$}
} & \parbox{32pt}{\raggedright
0.05106
} & \parbox{38pt}{\raggedright
\small{8.78$\times10^{-8}$}
} & \parbox{32pt}{\raggedright
0.05534
} & \parbox{38pt}{\raggedright
\small{4.03$\times10^{-9}$}
} & \parbox{32pt}{\raggedright
0.04229
} & \parbox{38pt}{\raggedright
\small{5.96$\times10^{-8}$}
} \\
\hline
\parbox{41pt}{\raggedright
(20,2)
} & \parbox{32pt}{\raggedright
0.05271
} & \parbox{38pt}{\raggedright
0.03781
} & \parbox{32pt}{\raggedright
0.10129
} & \parbox{38pt}{\raggedright
0.0117
} & \parbox{32pt}{\raggedright
0.03855
} & \parbox{38pt}{\raggedright
0.02808
} & \parbox{32pt}{\raggedright
0.07209
} & \parbox{38pt}{\raggedright
0.0061
} \\
\hline
\parbox{41pt}{\raggedright
(25,20)
} & \parbox{32pt}{\raggedright
0.06519
} & \parbox{38pt}{\raggedright
\small{2.07$\times10^{-9}$}
} & \parbox{32pt}{\raggedright
0.06472
} & \parbox{38pt}{\raggedright
\small{4.57$\times10^{-8}$}
} & \parbox{32pt}{\raggedright
0.05211
} & \parbox{38pt}{\raggedright
\small{2.07$\times10^{-9}$}
} & \parbox{32pt}{\raggedright
0.05178
} & \parbox{38pt}{\raggedright
\small{3.09$\times10^{-8}$}
} \\
\hline
\parbox{41pt}{\raggedright
(50,50)
} & \parbox{32pt}{\raggedright
0.0351
} & \parbox{38pt}{\raggedright
0
} & \parbox{32pt}{\raggedright
0.07716
} & \parbox{38pt}{\raggedright
0
} & \parbox{32pt}{\raggedright
0.03036
} & \parbox{38pt}{\raggedright
0
} & \parbox{32pt}{\raggedright
0.0594
} & \parbox{38pt}{\raggedright
0
} \\
\hline
\parbox{41pt}{\raggedright
(50,8)
} & \parbox{32pt}{\raggedright
0.07247
} & \parbox{38pt}{\raggedright
0.00033
} & \parbox{32pt}{\raggedright
0.16596
} & \parbox{38pt}{\raggedright
0.00108
} & \parbox{32pt}{\raggedright
0.05529
} & \parbox{38pt}{\raggedright
0.00031
} & \parbox{32pt}{\raggedright
0.11985
} & \parbox{38pt}{\raggedright
0.00066
} \\
\hline
\parbox{41pt}{\raggedright
(60,50)
} & \parbox{32pt}{\raggedright
0.03254
} & \parbox{38pt}{\raggedright
0
} & \parbox{32pt}{\raggedright
0.09167
} & \parbox{38pt}{\raggedright
0
} & \parbox{32pt}{\raggedright
0.02837
} & \parbox{38pt}{\raggedright
0
} & \parbox{32pt}{\raggedright
0.06919
} & \parbox{38pt}{\raggedright
0
} \\
\hline
\parbox{41pt}{\raggedright
(100,100)
} & \parbox{32pt}{\raggedright
0.01879
} & \parbox{38pt}{\raggedright
0
} & \parbox{32pt}{\raggedright
0.10239
} & \parbox{38pt}{\raggedright
0
} & \parbox{32pt}{\raggedright
0.01724
} & \parbox{38pt}{\raggedright
0
} & \parbox{32pt}{\raggedright
0.07495
} & \parbox{38pt}{\raggedright
0
} \\
\hline
\parbox{41pt}{\raggedright
(4,100)
} & \parbox{32pt}{\raggedright
0.06559
} & \parbox{38pt}{\raggedright
0.01236
} & \parbox{32pt}{\raggedright
0.00823
} & \parbox{38pt}{\raggedright
0.01024
} & \parbox{32pt}{\raggedright
0.04804
} & \parbox{38pt}{\raggedright
0.01003
} & \parbox{32pt}{\raggedright
0.00293
} & \parbox{38pt}{\raggedright
0.00557
} \\
\hline
\parbox{41pt}{\raggedright
(100,4)
} & \parbox{32pt}{\raggedright
0.06559
} & \parbox{38pt}{\raggedright
0.01236
} & \parbox{32pt}{\raggedright
0.26446
} & \parbox{38pt}{\raggedright
0.01024
} & \parbox{32pt}{\raggedright
0.04804
} & \parbox{38pt}{\raggedright
0.01003
} & \parbox{32pt}{\raggedright
0.18535
} & \parbox{38pt}{\raggedright
0.00557
} \\
\hline
\parbox{41pt}{\raggedright
(150,100)
} & \parbox{32pt}{\raggedright
0.01579
} & \parbox{38pt}{\raggedright
0
} & \parbox{32pt}{\raggedright
0.14059
} & \parbox{38pt}{\raggedright
0
} & \parbox{32pt}{\raggedright
0.01465
} & \parbox{38pt}{\raggedright
0
} & \parbox{32pt}{\raggedright
0.10082
} & \parbox{38pt}{\raggedright
0
} \\
\hline
\parbox{41pt}{\raggedright
(200,200)
} & \parbox{32pt}{\raggedright
0.00965
} & \parbox{38pt}{\raggedright
0
} & \parbox{32pt}{\raggedright
0.12808
} & \parbox{38pt}{\raggedright
0
} & \parbox{32pt}{\raggedright
0.00924
} & \parbox{38pt}{\raggedright
0
} & \parbox{32pt}{\raggedright
0.09101
} & \parbox{38pt}{\raggedright
0
} \\
\hline
\parbox{41pt}{\raggedright
(200,320)
} & \parbox{32pt}{\raggedright
0.00787
} & \parbox{38pt}{\raggedright
0
} & \parbox{32pt}{\raggedright
0.10243
} & \parbox{38pt}{\raggedright
0
} & \parbox{32pt}{\raggedright
0.00759
} & \parbox{38pt}{\raggedright
0
} & \parbox{32pt}{\raggedright
0.0729
} & \parbox{38pt}{\raggedright
0
} \\
\hline
\parbox{41pt}{\raggedright
(4,200)
} & \parbox{32pt}{\raggedright
0.06242
} & \parbox{38pt}{\raggedright
0.01236
} & \parbox{32pt}{\raggedright
0.00869
} & \parbox{38pt}{\raggedright
0.01024
} & \parbox{32pt}{\raggedright
0.0453
} & \parbox{38pt}{\raggedright
0.01003
} & \parbox{32pt}{\raggedright
0.00379
} & \parbox{38pt}{\raggedright
0.00557
} \\
\hline
\parbox{41pt}{\raggedright
(200,4)
} & \parbox{32pt}{\raggedright
0.06242
} & \parbox{38pt}{\raggedright
0.01236
} & \parbox{32pt}{\raggedright
0.35623
} & \parbox{38pt}{\raggedright
0.01024
} & \parbox{32pt}{\raggedright
0.0453
} & \parbox{38pt}{\raggedright
0.01003
} & \parbox{32pt}{\raggedright
0.24816
} & \parbox{38pt}{\raggedright
0.00557
} \\
\hline
\end{tabular}
\end{table}

\begin{table}\centering
\caption{\c{Absolute errors between actual values and their asymptotics for the case of $\rho=-1$}}
\label{table:2}

\vspace{0.5cm}
\begin{tabular}{|p{41pt}|p{32pt}|p{38pt}|p{32pt}|p{38pt}|p{32pt}|p{38pt}|p{32pt}|p{38pt}|}
\hline
\parbox{41pt}{\raggedright \multirow{2}{*}{
(x,y)
}} & \multicolumn{4}{|l|}{\parbox{140pt}{\raggedright
n=1000
}} & \multicolumn{4}{|l|}{\parbox{140pt}{\raggedright
n=10000
}} \\
\cline{2-9}
 & \parbox{32pt}{\raggedright
$\Delta_1^{p}$
} & \parbox{38pt}{\raggedright
$\Delta_1^{l}$
} & \parbox{32pt}{\raggedright
$\Delta_3^{p}$
} & \parbox{38pt}{\raggedright
$\Delta_3^{l}$
} & \parbox{32pt}{\raggedright
$\Delta_1^{p}$
} & \parbox{38pt}{\raggedright
$\Delta_1^{l}$
} & \parbox{32pt}{\raggedright
$\Delta_3^{p}$
} & \parbox{38pt}{\raggedright
$\Delta_3^{l}$
} \\
\hline
\parbox{41pt}{\raggedright
(0.5,0.5)
} & \parbox{32pt}{\raggedright
0.00112
} & \parbox{38pt}{\raggedright
0.02047
} & \parbox{32pt}{\raggedright
0.00051
} & \parbox{38pt}{\raggedright
0.00314
} & \parbox{32pt}{\raggedright
0.00079
} & \parbox{38pt}{\raggedright
0.01474
} & \parbox{32pt}{\raggedright
0.00033
} & \parbox{38pt}{\raggedright
0.00156
} \\
\hline
\parbox{41pt}{\raggedright
(1,1)
} & \parbox{32pt}{\raggedright
0.00027
} & \parbox{38pt}{\raggedright
0.04775
} & \parbox{32pt}{\raggedright
0.00027
} & \parbox{38pt}{\raggedright
0.00763
} & \parbox{32pt}{\raggedright
0.00003
} & \parbox{38pt}{\raggedright
0.0343
} & \parbox{32pt}{\raggedright
0.00003
} & \parbox{38pt}{\raggedright
0.00393
} \\
\hline
\parbox{41pt}{\raggedright
(1,0.5)
} & \parbox{32pt}{\raggedright
0.00156
} & \parbox{38pt}{\raggedright
0.03177
} & \parbox{32pt}{\raggedright
0.00065
} & \parbox{38pt}{\raggedright
0.00502
} & \parbox{32pt}{\raggedright
0.00109
} & \parbox{38pt}{\raggedright
0.02285
} & \parbox{32pt}{\raggedright
0.00044
} & \parbox{38pt}{\raggedright
0.00256
} \\
\hline
\parbox{41pt}{\raggedright
(2,1.5)
} & \parbox{32pt}{\raggedright
0.02903
} & \parbox{38pt}{\raggedright
0.0664
} & \parbox{32pt}{\raggedright
0.00249
} & \parbox{38pt}{\raggedright
0.01606
} & \parbox{32pt}{\raggedright
0.02051
} & \parbox{38pt}{\raggedright
0.04861
} & \parbox{32pt}{\raggedright
0.00125
} & \parbox{38pt}{\raggedright
0.00833
} \\
\hline
\parbox{41pt}{\raggedright
(3,3)
} & \parbox{32pt}{\raggedright
0.07729
} & \parbox{38pt}{\raggedright
0.04709
} & \parbox{32pt}{\raggedright
0.00535
} & \parbox{38pt}{\raggedright
0.02369
} & \parbox{32pt}{\raggedright
0.05425
} & \parbox{38pt}{\raggedright
0.03635
} & \parbox{32pt}{\raggedright
0.00281
} & \parbox{38pt}{\raggedright
0.01253
} \\
\hline
\parbox{41pt}{\raggedright
(3,5)
} & \parbox{32pt}{\raggedright
0.09003
} & \parbox{38pt}{\raggedright
0.02961
} & \parbox{32pt}{\raggedright
0.00878
} & \parbox{38pt}{\raggedright
0.01902
} & \parbox{32pt}{\raggedright
0.06375
} & \parbox{38pt}{\raggedright
0.02331
} & \parbox{32pt}{\raggedright
0.00448
} & \parbox{38pt}{\raggedright
0.01026
} \\
\hline
\parbox{41pt}{\raggedright
(7,3)
} & \parbox{32pt}{\raggedright
0.09138
} & \parbox{38pt}{\raggedright
0.02526
} & \parbox{32pt}{\raggedright
0.01187
} & \parbox{38pt}{\raggedright
0.01477
} & \parbox{32pt}{\raggedright
0.06532
} & \parbox{38pt}{\raggedright
0.01968
} & \parbox{32pt}{\raggedright
0.00597
} & \parbox{38pt}{\raggedright
0.00796
} \\
\hline
\parbox{41pt}{\raggedright
(4,5)
} & \parbox{32pt}{\raggedright
0.10037
} & \parbox{38pt}{\raggedright
0.01771
} & \parbox{32pt}{\raggedright
0.01094
} & \parbox{38pt}{\raggedright
0.01678
} & \parbox{32pt}{\raggedright
0.07132
} & \parbox{38pt}{\raggedright
0.01458
} & \parbox{32pt}{\raggedright
0.00552
} & \parbox{38pt}{\raggedright
0.00923
} \\
\hline
\parbox{41pt}{\raggedright
(5,5)
} & \parbox{32pt}{\raggedright
0.10482
} & \parbox{38pt}{\raggedright
0.01089
} & \parbox{32pt}{\raggedright
0.01309
} & \parbox{38pt}{\raggedright
0.01347
} & \parbox{32pt}{\raggedright
0.07485
} & \parbox{38pt}{\raggedright
0.00928
} & \parbox{32pt}{\raggedright
0.00656
} & \parbox{38pt}{\raggedright
0.00755
} \\
\hline
\parbox{41pt}{\raggedright
(5,8)
} & \parbox{32pt}{\raggedright
0.10578
} & \parbox{38pt}{\raggedright
0.00579
} & \parbox{32pt}{\raggedright
0.01833
} & \parbox{38pt}{\raggedright
0.00786
} & \parbox{32pt}{\raggedright
0.07659
} & \parbox{38pt}{\raggedright
0.00497
} & \parbox{32pt}{\raggedright
0.0091
} & \parbox{38pt}{\raggedright
0.00446
} \\
\hline
\parbox{41pt}{\raggedright
(6,7)
} & \parbox{32pt}{\raggedright
0.10798
} & \parbox{38pt}{\raggedright
0.0031
} & \parbox{32pt}{\raggedright
0.01901
} & \parbox{38pt}{\raggedright
0.00611
} & \parbox{32pt}{\raggedright
0.07826
} & \parbox{38pt}{\raggedright
0.00279
} & \parbox{32pt}{\raggedright
0.00942
} & \parbox{38pt}{\raggedright
0.00356
} \\
\hline
\parbox{41pt}{\raggedright
(7,4)
} & \parbox{32pt}{\raggedright
0.10197
} & \parbox{38pt}{\raggedright
0.01321
} & \parbox{32pt}{\raggedright
0.01439
} & \parbox{38pt}{\raggedright
0.01232
} & \parbox{32pt}{\raggedright
0.07315
} & \parbox{38pt}{\raggedright
0.01082
} & \parbox{32pt}{\raggedright
0.0072
} & \parbox{38pt}{\raggedright
0.00681
} \\
\hline
\parbox{41pt}{\raggedright
(8,9)
} & \parbox{32pt}{\raggedright
0.10507
} & \parbox{38pt}{\raggedright
0.00045
} & \parbox{32pt}{\raggedright
0.02567
} & \parbox{38pt}{\raggedright
0.00159
} & \parbox{32pt}{\raggedright
0.07757
} & \parbox{38pt}{\raggedright
0.00043
} & \parbox{32pt}{\raggedright
0.01269
} & \parbox{38pt}{\raggedright
0.00098
} \\
\hline
\parbox{41pt}{\raggedright
(10,10)
} & \parbox{32pt}{\raggedright
0.10075
} & \parbox{38pt}{\raggedright
0.00009
} & \parbox{32pt}{\raggedright
0.02964
} & \parbox{38pt}{\raggedright
0.00048
} & \parbox{32pt}{\raggedright
0.07534
} & \parbox{38pt}{\raggedright
0.00009
} & \parbox{32pt}{\raggedright
0.01468
} & \parbox{38pt}{\raggedright
0.0003
} \\
\hline
\parbox{41pt}{\raggedright
(10,20)
} & \parbox{32pt}{\raggedright
0.08699
} & \parbox{38pt}{\raggedright
0.00004
} & \parbox{32pt}{\raggedright
0.03549
} & \parbox{38pt}{\raggedright
0.00024
} & \parbox{32pt}{\raggedright
0.06674
} & \parbox{38pt}{\raggedright
0.00004
} & \parbox{32pt}{\raggedright
0.01783
} & \parbox{38pt}{\raggedright
0.00015
} \\
\hline
\parbox{41pt}{\raggedright
(7,10)
} & \parbox{32pt}{\raggedright
0.10458
} & \parbox{38pt}{\raggedright
0.00091
} & \parbox{32pt}{\raggedright
0.02515
} & \parbox{38pt}{\raggedright
0.00238
} & \parbox{32pt}{\raggedright
0.07712
} & \parbox{38pt}{\raggedright
0.00085
} & \parbox{32pt}{\raggedright
0.01245
} & \parbox{38pt}{\raggedright
0.00142
} \\
\hline
\parbox{41pt}{\raggedright
(20,20)
} & \parbox{32pt}{\raggedright
0.07196
} & \parbox{38pt}{\raggedright
\small{4.12$\times10^{-9}$}
} & \parbox{32pt}{\raggedright
0.04146
} & \parbox{38pt}{\raggedright
\small{9.08$\times10^{-8}$}
} & \parbox{32pt}{\raggedright
0.05723
} & \parbox{38pt}{\raggedright
\small{4.12$\times10^{-9}$}
} & \parbox{32pt}{\raggedright
0.02108
} & \parbox{38pt}{\raggedright
\small{6.14$\times10^{-8}$}
} \\
\hline
\parbox{41pt}{\raggedright
(20,2)
} & \parbox{32pt}{\raggedright
0.05603
} & \parbox{38pt}{\raggedright
0.03781
} & \parbox{32pt}{\raggedright
0.01558
} & \parbox{38pt}{\raggedright
0.0117
} & \parbox{32pt}{\raggedright
0.04149
} & \parbox{38pt}{\raggedright
0.02808
} & \parbox{32pt}{\raggedright
0.00795
} & \parbox{38pt}{\raggedright
0.0061
} \\
\hline
\parbox{41pt}{\raggedright
(25,20)
} & \parbox{32pt}{\raggedright
0.06702
} & \parbox{38pt}{\raggedright
\small{2.07$\times10^{-9}$}
} & \parbox{32pt}{\raggedright
0.04225
} & \parbox{38pt}{\raggedright
\small{4.59$\times10^{-8}$}
} & \parbox{32pt}{\raggedright
0.05384
} & \parbox{38pt}{\raggedright
\small{2.08$\times10^{-9}$}
} & \parbox{32pt}{\raggedright
0.0216
} & \parbox{38pt}{\raggedright
\small{3.1$\times10^{-8}$}
} \\
\hline
\parbox{41pt}{\raggedright
(30,30)
} & \parbox{32pt}{\raggedright
0.05432
} & \parbox{38pt}{\raggedright
\small{1.9$\times10^{-13}$}
} & \parbox{32pt}{\raggedright
0.04344
} & \parbox{38pt}{\raggedright
\small{9.2$\times10^{-12}$}
} & \parbox{32pt}{\raggedright
0.04492
} & \parbox{38pt}{\raggedright
\small{1.9$\times10^{-13}$}
} & \parbox{32pt}{\raggedright
0.02258
} & \parbox{38pt}{\raggedright
\small{6.3$\times10^{-12}$}
} \\
\hline
\parbox{41pt}{\raggedright
(35,40)
} & \parbox{32pt}{\raggedright
0.04581
} & \parbox{38pt}{\raggedright
\small{6.7$\times10^{-16}$}
} & \parbox{32pt}{\raggedright
0.04304
} & \parbox{38pt}{\raggedright
\small{4.2$\times10^{-14}$}
} & \parbox{32pt}{\raggedright
0.03868
} & \parbox{38pt}{\raggedright
\small{6.7$\times10^{-16}$}
} & \parbox{32pt}{\raggedright
0.02267
} & \parbox{38pt}{\raggedright
\small{2.9$\times10^{-14}$}
} \\
\hline
\parbox{41pt}{\raggedright
(40,40)
} & \parbox{32pt}{\raggedright
0.04335
} & \parbox{38pt}{\raggedright
0
} & \parbox{32pt}{\raggedright
0.04279
} & \parbox{38pt}{\raggedright
\small{6.7$\times10^{-16}$}
} & \parbox{32pt}{\raggedright
0.03685
} & \parbox{38pt}{\raggedright
0
} & \parbox{32pt}{\raggedright
0.02263
} & \parbox{38pt}{\raggedright
\small{4.4$\times10^{-16}$}
} \\
\hline
\parbox{41pt}{\raggedright
(50,50)
} & \parbox{32pt}{\raggedright
0.03597
} & \parbox{38pt}{\raggedright
0
} & \parbox{32pt}{\raggedright
0.04136
} & \parbox{38pt}{\raggedright
0
} & \parbox{32pt}{\raggedright
0.03122
} & \parbox{38pt}{\raggedright
0
} & \parbox{32pt}{\raggedright
0.02218
} & \parbox{38pt}{\raggedright
0
} \\
\hline
\parbox{41pt}{\raggedright
(55,60)
} & \parbox{32pt}{\raggedright
0.03191
} & \parbox{38pt}{\raggedright
0
} & \parbox{32pt}{\raggedright
0.04014
} & \parbox{38pt}{\raggedright
0
} & \parbox{32pt}{\raggedright
0.02804
} & \parbox{38pt}{\raggedright
0
} & \parbox{32pt}{\raggedright
0.02171
} & \parbox{38pt}{\raggedright
0
} \\
\hline
\parbox{41pt}{\raggedright
(150,100)
} & \parbox{32pt}{\raggedright
0.01616
} & \parbox{38pt}{\raggedright
0
} & \parbox{32pt}{\raggedright
0.0311
} & \parbox{38pt}{\raggedright
0
} & \parbox{32pt}{\raggedright
0.01501
} & \parbox{38pt}{\raggedright
0
} & \parbox{32pt}{\raggedright
0.01762
} & \parbox{38pt}{\raggedright
0
} \\
\hline
\parbox{41pt}{\raggedright
(200,200)
} & \parbox{32pt}{\raggedright
0.00988
} & \parbox{38pt}{\raggedright
0
} & \parbox{32pt}{\raggedright
0.02472
} & \parbox{38pt}{\raggedright
0
} & \parbox{32pt}{\raggedright
0.00946
} & \parbox{38pt}{\raggedright
0
} & \parbox{32pt}{\raggedright
0.01443
} & \parbox{38pt}{\raggedright
0
} \\
\hline
\end{tabular}
\end{table}

\begin{table}\centering
\caption{\c{Absolute errors between actual values and their asymptotics for the case of $\rho=0$}}
\label{table:3}

\vspace{0.5cm}
\begin{tabular}{|p{41pt}|p{32pt}|p{38pt}|p{32pt}|p{38pt}|p{32pt}|p{38pt}|p{32pt}|p{38pt}|}
\hline
\parbox{41pt}{\raggedright \multirow{2}{*}{
(x,y)
}} & \multicolumn{4}{|l|}{\parbox{140pt}{\raggedright
n=1000
}} & \multicolumn{4}{|l|}{\parbox{140pt}{\raggedright
n=10000
}} \\
\cline{2-9}
 & \parbox{32pt}{\raggedright
$\Delta_1^{p}$
} & \parbox{38pt}{\raggedright
$\Delta_1^{l}$
} & \parbox{32pt}{\raggedright
$\Delta_3^{p}$
} & \parbox{38pt}{\raggedright
$\Delta_3^{l}$
} & \parbox{32pt}{\raggedright
$\Delta_1^{p}$
} & \parbox{38pt}{\raggedright
$\Delta_1^{l}$
} & \parbox{32pt}{\raggedright
$\Delta_3^{p}$
} & \parbox{38pt}{\raggedright
$\Delta_3^{l}$
} \\
\hline
\parbox{41pt}{\raggedright
(0.5,0.5)
} & \parbox{32pt}{\raggedright
0.00105
} & \parbox{38pt}{\raggedright
0.02057
} & \parbox{32pt}{\raggedright
0.00059
} & \parbox{38pt}{\raggedright
0.00303
} & \parbox{32pt}{\raggedright
0.00079
} & \parbox{38pt}{\raggedright
0.01475
} & \parbox{32pt}{\raggedright
0.00034
} & \parbox{38pt}{\raggedright
0.00155
} \\
\hline
\parbox{41pt}{\raggedright
(1,1)
} & \parbox{32pt}{\raggedright
0.00014
} & \parbox{38pt}{\raggedright
0.0478
} & \parbox{32pt}{\raggedright
0.00014
} & \parbox{38pt}{\raggedright
0.00757
} & \parbox{32pt}{\raggedright
0.00001
} & \parbox{38pt}{\raggedright
0.03431
} & \parbox{32pt}{\raggedright
0.00001
} & \parbox{38pt}{\raggedright
0.00393
} \\
\hline
\parbox{41pt}{\raggedright
(1,0.5)
} & \parbox{32pt}{\raggedright
0.00147
} & \parbox{38pt}{\raggedright
0.03184
} & \parbox{32pt}{\raggedright
0.00075
} & \parbox{38pt}{\raggedright
0.00495
} & \parbox{32pt}{\raggedright
0.00108
} & \parbox{38pt}{\raggedright
0.02285
} & \parbox{32pt}{\raggedright
0.00045
} & \parbox{38pt}{\raggedright
0.00255
} \\
\hline
\parbox{41pt}{\raggedright
(2,1.5)
} & \parbox{32pt}{\raggedright
0.02913
} & \parbox{38pt}{\raggedright
0.06642
} & \parbox{32pt}{\raggedright
0.00239
} & \parbox{38pt}{\raggedright
0.01605
} & \parbox{32pt}{\raggedright
0.02052
} & \parbox{38pt}{\raggedright
0.04861
} & \parbox{32pt}{\raggedright
0.00124
} & \parbox{38pt}{\raggedright
0.00833
} \\
\hline
\parbox{41pt}{\raggedright
(3,3)
} & \parbox{32pt}{\raggedright
0.07733
} & \parbox{38pt}{\raggedright
0.04709
} & \parbox{32pt}{\raggedright
0.00531
} & \parbox{38pt}{\raggedright
0.02369
} & \parbox{32pt}{\raggedright
0.05425
} & \parbox{38pt}{\raggedright
0.03635
} & \parbox{32pt}{\raggedright
0.0028
} & \parbox{38pt}{\raggedright
0.01253
} \\
\hline
\parbox{41pt}{\raggedright
(3,5)
} & \parbox{32pt}{\raggedright
0.09005
} & \parbox{38pt}{\raggedright
0.02961
} & \parbox{32pt}{\raggedright
0.00876
} & \parbox{38pt}{\raggedright
0.01902
} & \parbox{32pt}{\raggedright
0.06375
} & \parbox{38pt}{\raggedright
0.02331
} & \parbox{32pt}{\raggedright
0.00447
} & \parbox{38pt}{\raggedright
0.01026
} \\
\hline
\parbox{41pt}{\raggedright
(2,3)
} & \parbox{32pt}{\raggedright
0.05766
} & \parbox{38pt}{\raggedright
0.05824
} & \parbox{32pt}{\raggedright
0.00402
} & \parbox{38pt}{\raggedright
0.02136
} & \parbox{32pt}{\raggedright
0.04046
} & \parbox{38pt}{\raggedright
0.04377
} & \parbox{32pt}{\raggedright
0.00213
} & \parbox{38pt}{\raggedright
0.0119
} \\
\hline
\parbox{41pt}{\raggedright
(4,5)
} & \parbox{32pt}{\raggedright
0.10038
} & \parbox{38pt}{\raggedright
0.01771
} & \parbox{32pt}{\raggedright
0.01092
} & \parbox{38pt}{\raggedright
0.01678
} & \parbox{32pt}{\raggedright
0.07133
} & \parbox{38pt}{\raggedright
0.01458
} & \parbox{32pt}{\raggedright
0.00552
} & \parbox{38pt}{\raggedright
0.00923
} \\
\hline
\parbox{41pt}{\raggedright
(5,5)
} & \parbox{32pt}{\raggedright
0.10483
} & \parbox{38pt}{\raggedright
0.01089
} & \parbox{32pt}{\raggedright
0.01309
} & \parbox{38pt}{\raggedright
0.01347
} & \parbox{32pt}{\raggedright
0.07485
} & \parbox{38pt}{\raggedright
0.00928
} & \parbox{32pt}{\raggedright
0.00656
} & \parbox{38pt}{\raggedright
0.00755
} \\
\hline
\parbox{41pt}{\raggedright
(5,8)
} & \parbox{32pt}{\raggedright
0.10579

} & \parbox{38pt}{\raggedright
0.00579
} & \parbox{32pt}{\raggedright
0.01833
} & \parbox{38pt}{\raggedright
0.00786
} & \parbox{32pt}{\raggedright
0.07659
} & \parbox{38pt}{\raggedright
0.00497
} & \parbox{32pt}{\raggedright
0.0091
} & \parbox{38pt}{\raggedright
0.00446
} \\
\hline
\parbox{41pt}{\raggedright
(6,7)
} & \parbox{32pt}{\raggedright
0.10799
} & \parbox{38pt}{\raggedright
0.0031
} & \parbox{32pt}{\raggedright
0.019
} & \parbox{38pt}{\raggedright
0.00611
} & \parbox{32pt}{\raggedright
0.07826
} & \parbox{38pt}{\raggedright
0.00279
} & \parbox{32pt}{\raggedright
0.00942
} & \parbox{38pt}{\raggedright
0.00356
} \\
\hline
\parbox{41pt}{\raggedright
(7,4)
} & \parbox{32pt}{\raggedright
0.10198
} & \parbox{38pt}{\raggedright
0.01321
} & \parbox{32pt}{\raggedright
0.01438
} & \parbox{38pt}{\raggedright
0.01232
} & \parbox{32pt}{\raggedright
0.07315
} & \parbox{38pt}{\raggedright
0.01082
} & \parbox{32pt}{\raggedright
0.0072
} & \parbox{38pt}{\raggedright
0.00681
} \\
\hline
\parbox{41pt}{\raggedright
(8,9)
} & \parbox{32pt}{\raggedright
0.10507
} & \parbox{38pt}{\raggedright
0.00045
} & \parbox{32pt}{\raggedright
0.02567
} & \parbox{38pt}{\raggedright
0.00159
} & \parbox{32pt}{\raggedright
0.07757
} & \parbox{38pt}{\raggedright
0.00043
} & \parbox{32pt}{\raggedright
0.01269
} & \parbox{38pt}{\raggedright
0.00098
} \\
\hline
\parbox{41pt}{\raggedright
(10,10)
} & \parbox{32pt}{\raggedright
0.10075
} & \parbox{38pt}{\raggedright
0.00009
} & \parbox{32pt}{\raggedright
0.02964
} & \parbox{38pt}{\raggedright
0.00048
} & \parbox{32pt}{\raggedright
0.07534
} & \parbox{38pt}{\raggedright
0.00009
} & \parbox{32pt}{\raggedright
0.01468
} & \parbox{38pt}{\raggedright
0.0003
} \\
\hline
\parbox{41pt}{\raggedright
(10,20)
} & \parbox{32pt}{\raggedright
0.087
} & \parbox{38pt}{\raggedright
0.00004
} & \parbox{32pt}{\raggedright
0.03549
} & \parbox{38pt}{\raggedright
0.00024
} & \parbox{32pt}{\raggedright
0.06674
} & \parbox{38pt}{\raggedright
0.00004
} & \parbox{32pt}{\raggedright
0.01783
} & \parbox{38pt}{\raggedright
0.00015
} \\
\hline
\parbox{41pt}{\raggedright
(7,10)
} & \parbox{32pt}{\raggedright
0.10458
} & \parbox{38pt}{\raggedright
0.00091
} & \parbox{32pt}{\raggedright
0.02515
} & \parbox{38pt}{\raggedright
0.00238
} & \parbox{32pt}{\raggedright
0.07712
} & \parbox{38pt}{\raggedright
0.00085
} & \parbox{32pt}{\raggedright
0.01245
} & \parbox{38pt}{\raggedright
0.00142
} \\
\hline
\parbox{41pt}{\raggedright
(20,20)
} & \parbox{32pt}{\raggedright
0.07196
} & \parbox{38pt}{\raggedright
\small{4.12$\times10^{-9}$}
} & \parbox{32pt}{\raggedright
0.04146
} & \parbox{38pt}{\raggedright
\small{9.08$\times10^{-8}$}
} & \parbox{32pt}{\raggedright
0.05723
} & \parbox{38pt}{\raggedright
\small{4.12$\times10^{-9}$}
} & \parbox{32pt}{\raggedright
0.02108
} & \parbox{38pt}{\raggedright
\small{6.14$\times10^{-8}$}
} \\
\hline
\parbox{41pt}{\raggedright
(20,2)
} & \parbox{32pt}{\raggedright
0.05603
} & \parbox{38pt}{\raggedright
0.03781
} & \parbox{32pt}{\raggedright
0.01558
} & \parbox{38pt}{\raggedright
0.0117
} & \parbox{32pt}{\raggedright
0.04149
} & \parbox{38pt}{\raggedright
0.02808
} & \parbox{32pt}{\raggedright
0.00795
} & \parbox{38pt}{\raggedright
0.0061
} \\
\hline
\parbox{41pt}{\raggedright
(25,20)
} & \parbox{32pt}{\raggedright
0.06702
} & \parbox{38pt}{\raggedright
\small{2.08$\times10^{-9}$}
} & \parbox{32pt}{\raggedright
0.04225
} & \parbox{38pt}{\raggedright
\small{4.59$\times10^{-8}$}
} & \parbox{32pt}{\raggedright
0.05384
} & \parbox{38pt}{\raggedright
\small{2.08$\times10^{-9}$}
} & \parbox{32pt}{\raggedright
0.0216
} & \parbox{38pt}{\raggedright
\small{3.1$\times10^{-8}$}
} \\
\hline
\parbox{41pt}{\raggedright
(30,30)
} & \parbox{32pt}{\raggedright
0.05432
} & \parbox{38pt}{\raggedright
\small{1.9$\times10^{-13}$}
} & \parbox{32pt}{\raggedright
0.04344
} & \parbox{38pt}{\raggedright
\small{9.2$\times10^{-12}$}
} & \parbox{32pt}{\raggedright
0.04492
} & \parbox{38pt}{\raggedright
\small{1.9$\times10^{-13}$}
} & \parbox{32pt}{\raggedright
0.02258
} & \parbox{38pt}{\raggedright
\small{6.3$\times10^{-12}$}
} \\
\hline
\parbox{41pt}{\raggedright
(35,40)
} & \parbox{32pt}{\raggedright
0.04581
} & \parbox{38pt}{\raggedright
\small{6.7$\times10^{-16}$}
} & \parbox{32pt}{\raggedright
0.04304
} & \parbox{38pt}{\raggedright
\small{4.2$\times10^{-14}$}
} & \parbox{32pt}{\raggedright
0.03868
} & \parbox{38pt}{\raggedright
\small{6.7$\times10^{-16}$}
} & \parbox{32pt}{\raggedright
0.02267
} & \parbox{38pt}{\raggedright
\small{2.9$\times10^{-14}$}
} \\
\hline
\parbox{41pt}{\raggedright
(40,40)
} & \parbox{32pt}{\raggedright
0.04335
} & \parbox{38pt}{\raggedright
0
} & \parbox{32pt}{\raggedright
0.04279
} & \parbox{38pt}{\raggedright
\small{6.7$\times10^{-16}$}
} & \parbox{32pt}{\raggedright
0.03685
} & \parbox{38pt}{\raggedright
0
} & \parbox{32pt}{\raggedright
0.02263
} & \parbox{38pt}{\raggedright
\small{4.4$\times10^{-16}$}
} \\
\hline
\parbox{41pt}{\raggedright
(50,50)
} & \parbox{32pt}{\raggedright
0.03597
} & \parbox{38pt}{\raggedright
0
} & \parbox{32pt}{\raggedright
0.04136
} & \parbox{38pt}{\raggedright
0
} & \parbox{32pt}{\raggedright
0.03122
} & \parbox{38pt}{\raggedright
0
} & \parbox{32pt}{\raggedright
0.02218
} & \parbox{38pt}{\raggedright
0
} \\
\hline
\parbox{41pt}{\raggedright
(55,60)
} & \parbox{32pt}{\raggedright
0.03191
} & \parbox{38pt}{\raggedright
0
} & \parbox{32pt}{\raggedright
0.04014
} & \parbox{38pt}{\raggedright
0
} & \parbox{32pt}{\raggedright
0.02804
} & \parbox{38pt}{\raggedright
0
} & \parbox{32pt}{\raggedright
0.02171
} & \parbox{38pt}{\raggedright
0
} \\
\hline
\parbox{41pt}{\raggedright
(150,100)
} & \parbox{32pt}{\raggedright
0.01616
} & \parbox{38pt}{\raggedright
0
} & \parbox{32pt}{\raggedright
0.0311
} & \parbox{38pt}{\raggedright
0
} & \parbox{32pt}{\raggedright
0.01501
} & \parbox{38pt}{\raggedright
0
} & \parbox{32pt}{\raggedright
0.01762
} & \parbox{38pt}{\raggedright
0
} \\
\hline
\parbox{41pt}{\raggedright
(200,200)
} & \parbox{32pt}{\raggedright
0.00988
} & \parbox{38pt}{\raggedright
0
} & \parbox{32pt}{\raggedright
0.02472
} & \parbox{38pt}{\raggedright
0
} & \parbox{32pt}{\raggedright
0.00946
} & \parbox{38pt}{\raggedright
0
} & \parbox{32pt}{\raggedright
0.01443
} & \parbox{38pt}{\raggedright
0
} \\
\hline
\end{tabular}
\end{table}

\begin{table}\centering
\caption{\c{Absolute errors between actual values and their asymptotics for the case of $\rho=1$}}
\label{table:4}

\vspace{0.5cm}
\begin{tabular}{|p{41pt}|p{32pt}|p{38pt}|p{32pt}|p{38pt}|p{32pt}|p{38pt}|p{32pt}|p{38pt}|}
\hline
\parbox{41pt}{\raggedright \multirow{2}{*}{
(x,y)
}} & \multicolumn{4}{|l|}{\parbox{140pt}{\raggedright
n=1000
}} & \multicolumn{4}{|l|}{\parbox{140pt}{\raggedright
n=10000
}} \\
\cline{2-9}
 & \parbox{32pt}{\raggedright
$\Delta_1^{p}$
} & \parbox{38pt}{\raggedright
$\Delta_1^{l}$
} & \parbox{32pt}{\raggedright
$\Delta_4^{p}$
} & \parbox{38pt}{\raggedright
$\Delta_4^{l}$
} & \parbox{32pt}{\raggedright
$\Delta_1^{p}$
} & \parbox{38pt}{\raggedright
$\Delta_1^{l}$
} & \parbox{32pt}{\raggedright
$\Delta_4^{p}$
} & \parbox{38pt}{\raggedright
$\Delta_4^{l}$
} \\
\hline
\parbox{41pt}{\raggedright
(0.5,0.5)
} & \parbox{32pt}{\raggedright
0.00392
} & \parbox{38pt}{\raggedright
0.01855
} & \parbox{32pt}{\raggedright
0.00211
} & \parbox{38pt}{\raggedright
0.0031
} & \parbox{32pt}{\raggedright
0.00294
} & \parbox{38pt}{\raggedright
0.01336
} & \parbox{32pt}{\raggedright
0.00123
} & \parbox{38pt}{\raggedright
0.00158
} \\
\hline
\parbox{41pt}{\raggedright
(1,1)
} & \parbox{32pt}{\raggedright
0.00018
} & \parbox{38pt}{\raggedright
0.03371
} & \parbox{32pt}{\raggedright
0.00018
} & \parbox{38pt}{\raggedright
0.00629
} & \parbox{32pt}{\raggedright
0.00002
} & \parbox{38pt}{\raggedright
0.02435
} & \parbox{32pt}{\raggedright
0.00002
} & \parbox{38pt}{\raggedright
0.00326
} \\
\hline
\parbox{41pt}{\raggedright
(1,0.5)
} & \parbox{32pt}{\raggedright
0.00392
} & \parbox{38pt}{\raggedright
0.01855
} & \parbox{32pt}{\raggedright
0.00211
} & \parbox{38pt}{\raggedright
0.0031
} & \parbox{32pt}{\raggedright
0.00294
} & \parbox{38pt}{\raggedright
0.01336
} & \parbox{32pt}{\raggedright
0.00123
} & \parbox{38pt}{\raggedright
0.00158
} \\
\hline
\parbox{41pt}{\raggedright
(2,1.5)
} & \parbox{32pt}{\raggedright
0.01828
} & \parbox{38pt}{\raggedright
0.03969
} & \parbox{32pt}{\raggedright
0.00214
} & \parbox{38pt}{\raggedright
0.00938
} & \parbox{32pt}{\raggedright
0.01301
} & \parbox{38pt}{\raggedright
0.02901
} & \parbox{32pt}{\raggedright
0.00109
} & \parbox{38pt}{\raggedright
0.00487
} \\
\hline
\parbox{41pt}{\raggedright
(3,3)
} & \parbox{32pt}{\raggedright
0.05207
} & \parbox{38pt}{\raggedright
0.02444
} & \parbox{32pt}{\raggedright
0.0056
} & \parbox{38pt}{\raggedright
0.01277
} & \parbox{32pt}{\raggedright
0.03691
} & \parbox{38pt}{\raggedright
0.01892
} & \parbox{32pt}{\raggedright
0.00291
} & \parbox{38pt}{\raggedright
0.00677
} \\
\hline
\parbox{41pt}{\raggedright
(3,5)
} & \parbox{32pt}{\raggedright
0.05207
} & \parbox{38pt}{\raggedright
0.02444
} & \parbox{32pt}{\raggedright
0.0056
} & \parbox{38pt}{\raggedright
0.01277
} & \parbox{32pt}{\raggedright
0.03691
} & \parbox{38pt}{\raggedright
0.01892
} & \parbox{32pt}{\raggedright
0.00291
} & \parbox{38pt}{\raggedright
0.00677
} \\
\hline
\parbox{41pt}{\raggedright
(2,3)
} & \parbox{32pt}{\raggedright
0.03393
} & \parbox{38pt}{\raggedright
0.03781
} & \parbox{32pt}{\raggedright
0.00334
} & \parbox{38pt}{\raggedright
0.0117
} & \parbox{32pt}{\raggedright
0.024
} & \parbox{38pt}{\raggedright
0.02808
} & \parbox{32pt}{\raggedright
0.00175
} & \parbox{38pt}{\raggedright
0.0061
} \\
\hline
\parbox{41pt}{\raggedright
(4,5)
} & \parbox{32pt}{\raggedright
0.05944
} & \parbox{38pt}{\raggedright
0.01236
} & \parbox{32pt}{\raggedright
0.008
} & \parbox{38pt}{\raggedright
0.01024
} & \parbox{32pt}{\raggedright
0.04248
} & \parbox{38pt}{\raggedright
0.01003
} & \parbox{32pt}{\raggedright
0.00409
} & \parbox{38pt}{\raggedright
0.00557
} \\
\hline
\parbox{41pt}{\raggedright
(5,5)
} & \parbox{32pt}{\raggedright
0.0617
} & \parbox{38pt}{\raggedright
0.00547
} & \parbox{32pt}{\raggedright
0.01032
} & \parbox{38pt}{\raggedright
0.0068
} & \parbox{32pt}{\raggedright
0.0445
} & \parbox{38pt}{\raggedright
0.00466
} & \parbox{32pt}{\raggedright
0.00522
} & \parbox{38pt}{\raggedright
0.00381
} \\
\hline
\parbox{41pt}{\raggedright
(5,8)
} & \parbox{32pt}{\raggedright
0.0617
} & \parbox{38pt}{\raggedright
0.00547
} & \parbox{32pt}{\raggedright
0.01032
} & \parbox{38pt}{\raggedright
0.0068
} & \parbox{32pt}{\raggedright
0.0445
} & \parbox{38pt}{\raggedright
0.00466
} & \parbox{32pt}{\raggedright
0.00522
} & \parbox{38pt}{\raggedright
0.00381
} \\
\hline
\parbox{41pt}{\raggedright
(6,7)
} & \parbox{32pt}{\raggedright
0.06152
} & \parbox{38pt}{\raggedright
0.00223
} & \parbox{32pt}{\raggedright
0.01238
} & \parbox{38pt}{\raggedright
0.00398
} & \parbox{32pt}{\raggedright
0.0448
} & \parbox{38pt}{\raggedright
0.00199
} & \parbox{32pt}{\raggedright
0.00622
} & \parbox{38pt}{\raggedright
0.0023
} \\
\hline
\parbox{41pt}{\raggedright
(8,7)
} & \parbox{32pt}{\raggedright
0.06019
} & \parbox{38pt}{\raggedright
0.00087
} & \parbox{32pt}{\raggedright
0.01415
} & \parbox{38pt}{\raggedright
0.00214
} & \parbox{32pt}{\raggedright
0.04423
} & \parbox{38pt}{\raggedright
0.00081
} & \parbox{32pt}{\raggedright
0.0071
} & \parbox{38pt}{\raggedright
0.00127
} \\
\hline
\parbox{41pt}{\raggedright
(8,9)
} & \parbox{32pt}{\raggedright
0.05831
} & \parbox{38pt}{\raggedright
0.00033
} & \parbox{32pt}{\raggedright
0.01566
} & \parbox{38pt}{\raggedright
0.00108
} & \parbox{32pt}{\raggedright
0.04323
} & \parbox{38pt}{\raggedright
0.00031
} & \parbox{32pt}{\raggedright
0.00785
} & \parbox{38pt}{\raggedright
0.00066
} \\
\hline
\parbox{41pt}{\raggedright
(10,10)
} & \parbox{32pt}{\raggedright
0.05406
} & \parbox{38pt}{\raggedright
0.00005
} & \parbox{32pt}{\raggedright
0.01799
} & \parbox{38pt}{\raggedright
0.00024
} & \parbox{32pt}{\raggedright
0.04072
} & \parbox{38pt}{\raggedright
0.00004
} & \parbox{32pt}{\raggedright
0.00903
} & \parbox{38pt}{\raggedright
0.00015
} \\
\hline
\parbox{41pt}{\raggedright
(10,20)
} & \parbox{32pt}{\raggedright
0.05406
} & \parbox{38pt}{\raggedright
0.00005
} & \parbox{32pt}{\raggedright
0.01799
} & \parbox{38pt}{\raggedright
0.00024
} & \parbox{32pt}{\raggedright
0.04072
} & \parbox{38pt}{\raggedright
0.00004
} & \parbox{32pt}{\raggedright
0.00903
} & \parbox{38pt}{\raggedright
0.00015
} \\
\hline
\parbox{41pt}{\raggedright
(7,10)
} & \parbox{32pt}{\raggedright
0.06019
} & \parbox{38pt}{\raggedright
0.00087
} & \parbox{32pt}{\raggedright
0.01415
} & \parbox{38pt}{\raggedright
0.00214
} & \parbox{32pt}{\raggedright
0.04423
} & \parbox{38pt}{\raggedright
0.00081
} & \parbox{32pt}{\raggedright
0.0071
} & \parbox{38pt}{\raggedright
0.00127
} \\
\hline
\parbox{41pt}{\raggedright
(20,20)
} & \parbox{32pt}{\raggedright
0.0371
} & \parbox{38pt}{\raggedright
\small{2.06$\times10^{-9}$}
} & \parbox{32pt}{\raggedright
0.02252
} & \parbox{38pt}{\raggedright
\small{4.54$\times10^{-8}$}
} & \parbox{32pt}{\raggedright
0.02962
} & \parbox{38pt}{\raggedright
\small{2.06$\times10^{-9}$}
} & \parbox{32pt}{\raggedright
0.01154
} & \parbox{38pt}{\raggedright
\small{3.07$\times10^{-8}$}
} \\
\hline
\parbox{41pt}{\raggedright
(20,2)
} & \parbox{32pt}{\raggedright
0.03393
} & \parbox{38pt}{\raggedright
0.03781
} & \parbox{32pt}{\raggedright
0.00334
} & \parbox{38pt}{\raggedright
0.0117
} & \parbox{32pt}{\raggedright
0.024
} & \parbox{38pt}{\raggedright
0.02808
} & \parbox{32pt}{\raggedright
0.00175
} & \parbox{38pt}{\raggedright
0.0061
} \\
\hline
\parbox{41pt}{\raggedright
(25,20)
} & \parbox{32pt}{\raggedright
0.0371
} & \parbox{38pt}{\raggedright
\small{2.06$\times10^{-9}$}
} & \parbox{32pt}{\raggedright
0.02252
} & \parbox{38pt}{\raggedright
\small{4.54$\times10^{-8}$}
} & \parbox{32pt}{\raggedright
0.02962
} & \parbox{38pt}{\raggedright
\small{2.06$\times10^{-9}$}
} & \parbox{32pt}{\raggedright
0.01539
} & \parbox{38pt}{\raggedright
\small{3.1$\times10^{-8}$}
} \\
\hline
\parbox{41pt}{\raggedright
(30,30)
} & \parbox{32pt}{\raggedright
0.02768
} & \parbox{38pt}{\raggedright
\small{9.4$\times10^{-14}$}
} & \parbox{32pt}{\raggedright
0.023
} & \parbox{38pt}{\raggedright
\small{4.6$\times10^{-12}$}
} & \parbox{32pt}{\raggedright
0.02295
} & \parbox{38pt}{\raggedright
\small{9.4$\times10^{-14}$}
} & \parbox{32pt}{\raggedright
0.01194
} & \parbox{38pt}{\raggedright
\small{3.2$\times10^{-12}$}
} \\
\hline
\parbox{41pt}{\raggedright
(35,40)
} & \parbox{32pt}{\raggedright
0.02451
} & \parbox{38pt}{\raggedright
\small{6.7$\times10^{-16}$}
} & \parbox{32pt}{\raggedright
0.02258
} & \parbox{38pt}{\raggedright
\small{4.2$\times10^{-14}$}
} & \parbox{32pt}{\raggedright
0.02062
} & \parbox{38pt}{\raggedright
\small{6.7$\times10^{-16}$}
} & \parbox{32pt}{\raggedright
0.0119
} & \parbox{38pt}{\raggedright
\small{2.9$\times10^{-14}$}
} \\
\hline
\parbox{41pt}{\raggedright
(40,40)
} & \parbox{32pt}{\raggedright
0.02198
} & \parbox{38pt}{\raggedright
0
} & \parbox{32pt}{\raggedright
0.02219
} & \parbox{38pt}{\raggedright
\small{4.4$\times10^{-16}$}
} & \parbox{32pt}{\raggedright
0.01871
} & \parbox{38pt}{\raggedright
0
} & \parbox{32pt}{\raggedright
0.01178
} & \parbox{38pt}{\raggedright
\small{2.2$\times10^{-16}$}
} \\
\hline
\parbox{41pt}{\raggedright
(50,50)
} & \parbox{32pt}{\raggedright
0.01818
} & \parbox{38pt}{\raggedright
0
} & \parbox{32pt}{\raggedright
0.02127
} & \parbox{38pt}{\raggedright
0
} & \parbox{32pt}{\raggedright
0.0158
} & \parbox{38pt}{\raggedright
0
} & \parbox{32pt}{\raggedright
0.01144
} & \parbox{38pt}{\raggedright
0
} \\
\hline
\parbox{41pt}{\raggedright
(55,60)
} & \parbox{32pt}{\raggedright
0.01673
} & \parbox{38pt}{\raggedright
0
} & \parbox{32pt}{\raggedright
0.020789
} & \parbox{38pt}{\raggedright
0
} & \parbox{32pt}{\raggedright
0.01466
} & \parbox{38pt}{\raggedright
0
} & \parbox{32pt}{\raggedright
0.01125
} & \parbox{38pt}{\raggedright
0
} \\
\hline
\parbox{41pt}{\raggedright
(150,100)
} & \parbox{32pt}{\raggedright
0.00967
} & \parbox{38pt}{\raggedright
0
} & \parbox{32pt}{\raggedright
0.01709
} & \parbox{38pt}{\raggedright
0
} & \parbox{32pt}{\raggedright
0.00889
} & \parbox{38pt}{\raggedright
0
} & \parbox{32pt}{\raggedright
0.00959
} & \parbox{38pt}{\raggedright
0
} \\
\hline
\parbox{41pt}{\raggedright
(200,200)
} & \parbox{32pt}{\raggedright
0.00495
} & \parbox{38pt}{\raggedright
0
} & \parbox{32pt}{\raggedright
0.01243
} & \parbox{38pt}{\raggedright
0
} & \parbox{32pt}{\raggedright
0.00474
} & \parbox{38pt}{\raggedright
0
} & \parbox{32pt}{\raggedright
0.00726
} & \parbox{38pt}{\raggedright
0
} \\
\hline
\end{tabular}
\end{table}

\begin{figure}\label{fig1}
\begin{center}
\subfigure[$\rho_n= -1$]
{%
 \epsfig{file=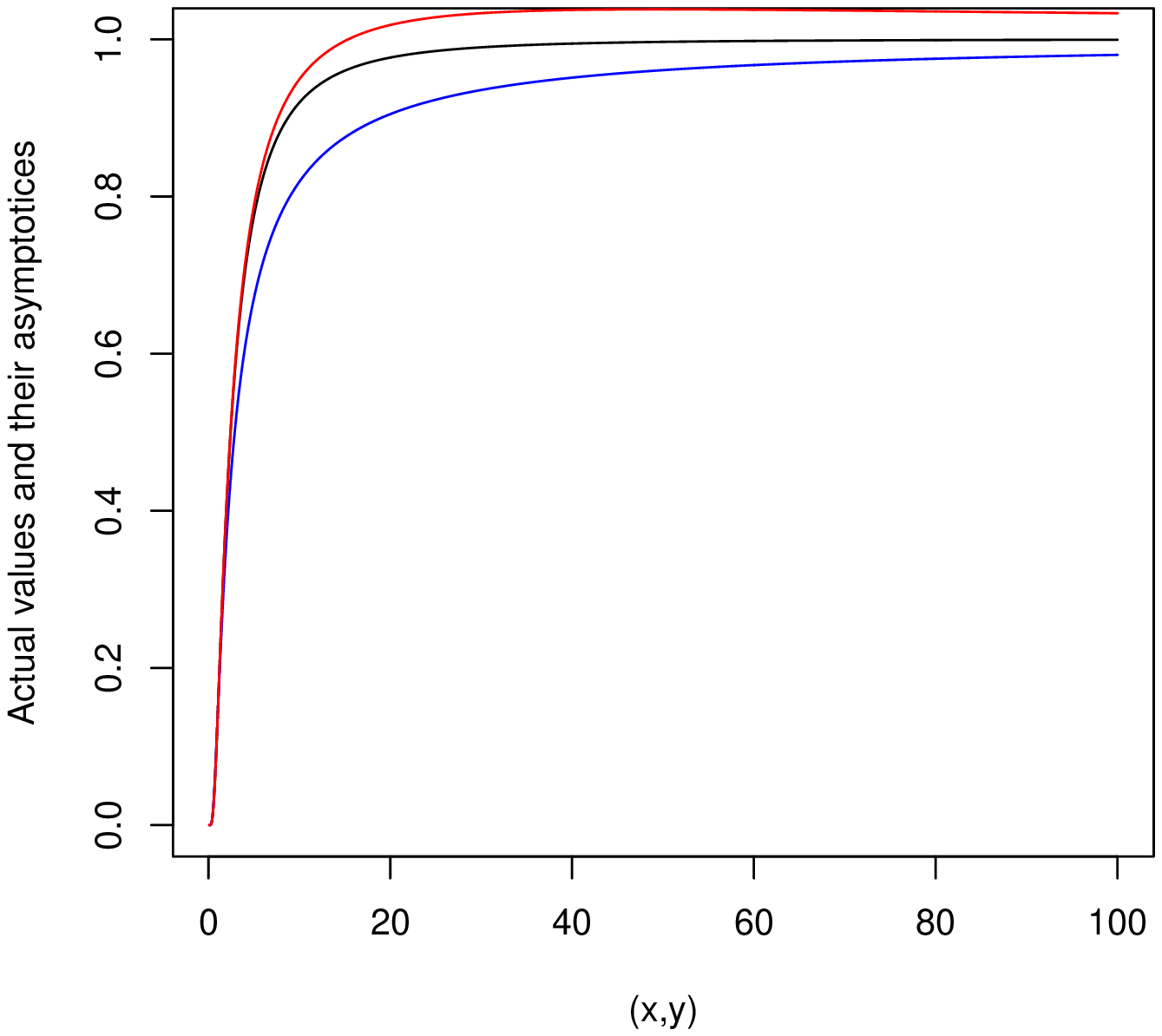, height=120pt, width=200pt,angle=0}
                   }%
\subfigure[$\lambda=2.5,\tau=-5, \rho_n=-0.915$]
{%
 \epsfig{file=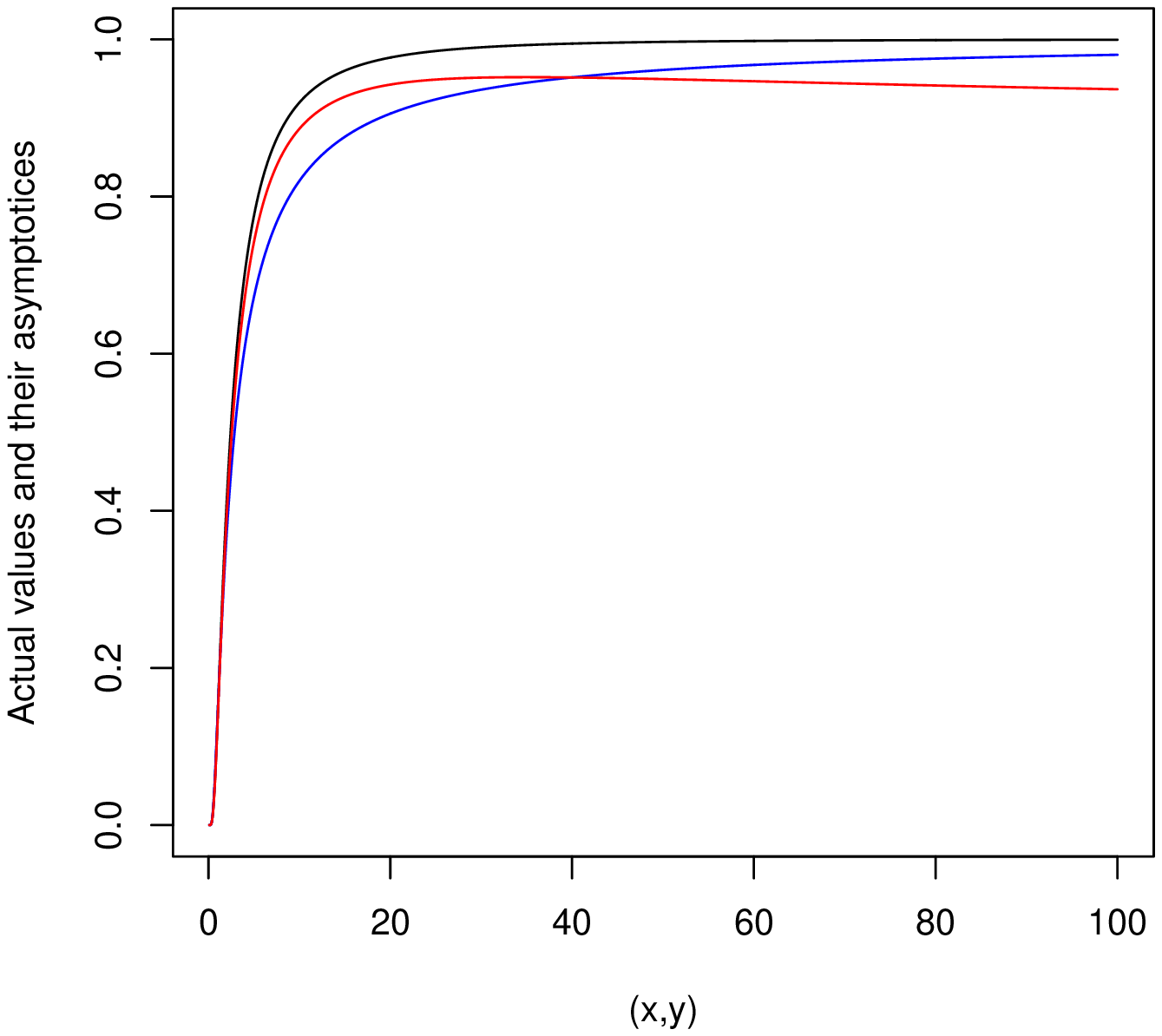, height=120pt, width=200pt,angle=0}
                   }%
                   \\
\subfigure[$\lambda=2.5,\tau=-2, \rho_n=-0.537$]
{%
 \epsfig{file=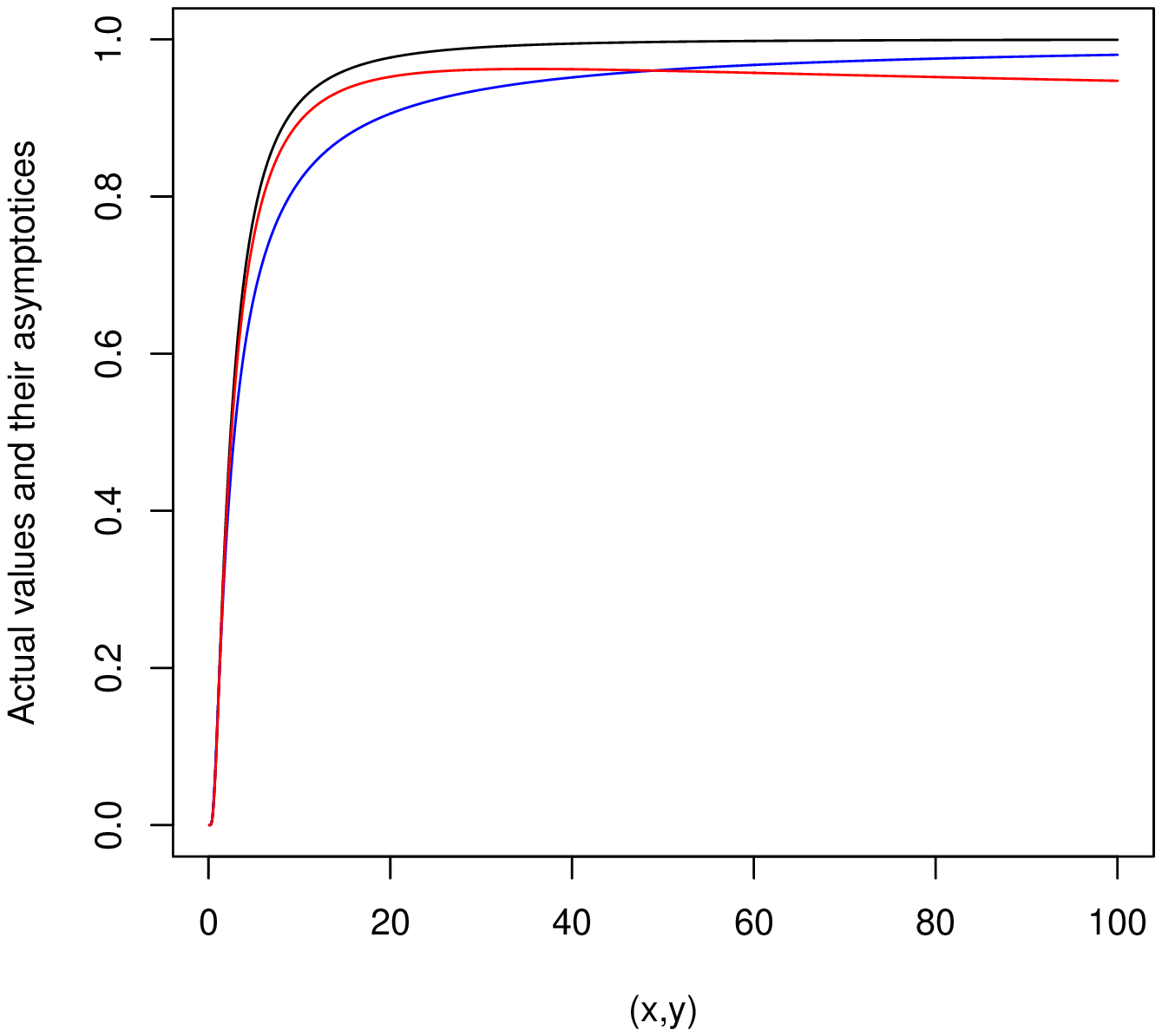, height=120pt, width=200pt,angle=0}
                   }%
 \subfigure[$\rho_n= 0$]
{%
 \epsfig{file=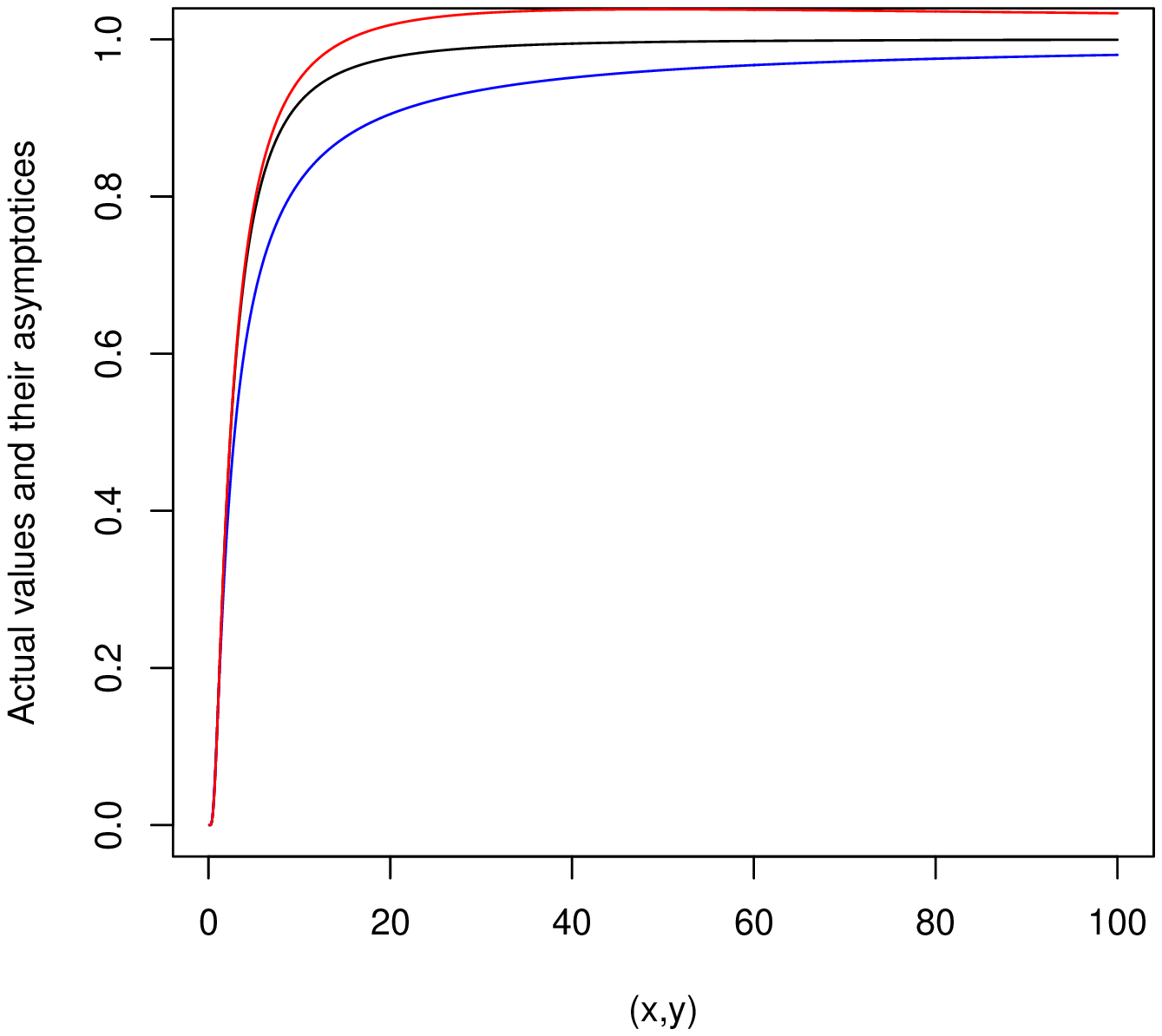, height=120pt, width=200pt,angle=0}
                   }%
                   \\
 \subfigure[$\lambda=2,\tau=2, \rho_n=0.329$]
{%
 \epsfig{file=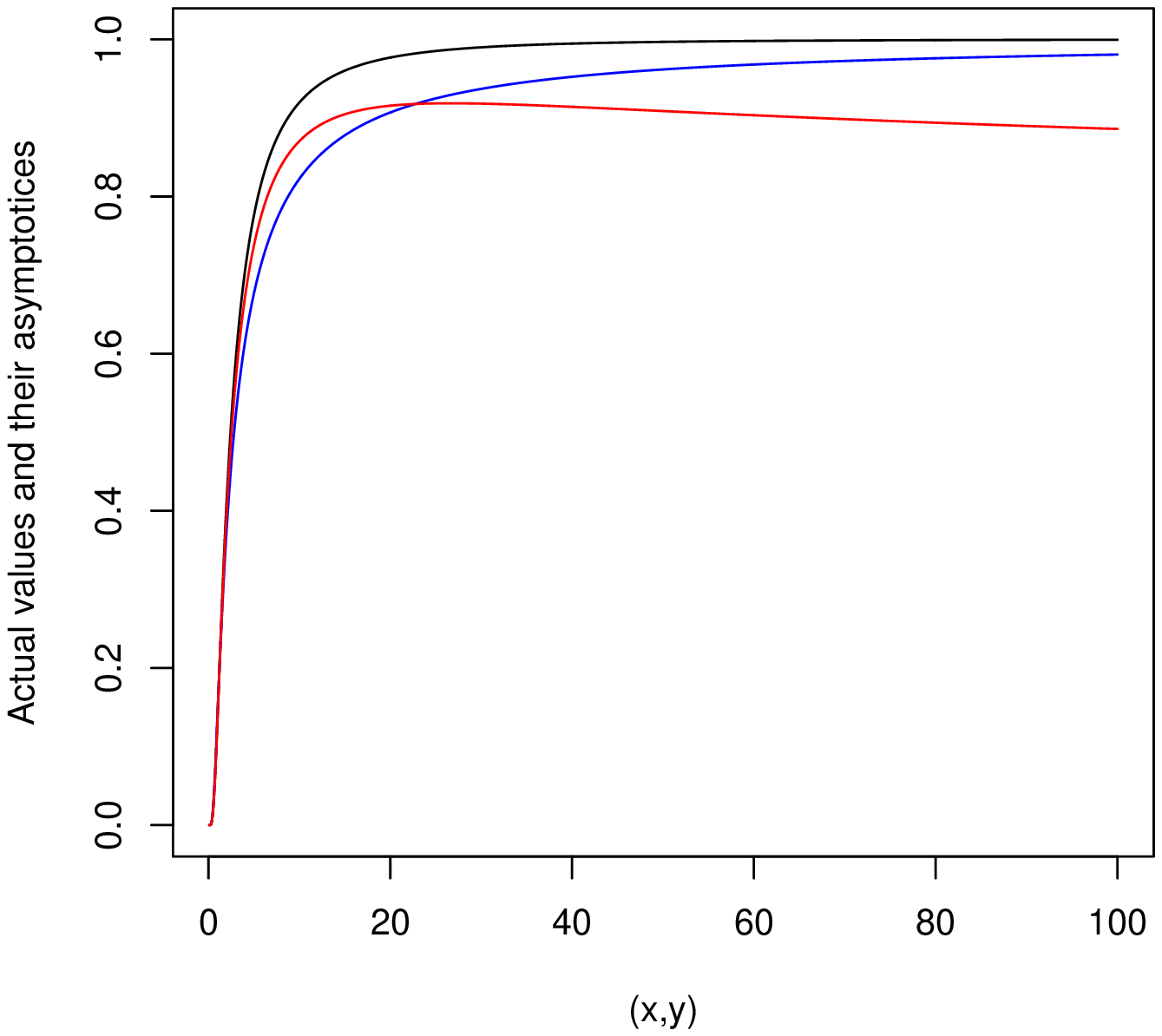, height=120pt, width=200pt,angle=0}
                   }%
\subfigure[$ \lambda=2,\tau=3,\rho_n=0.405$]
{%
\epsfig{file=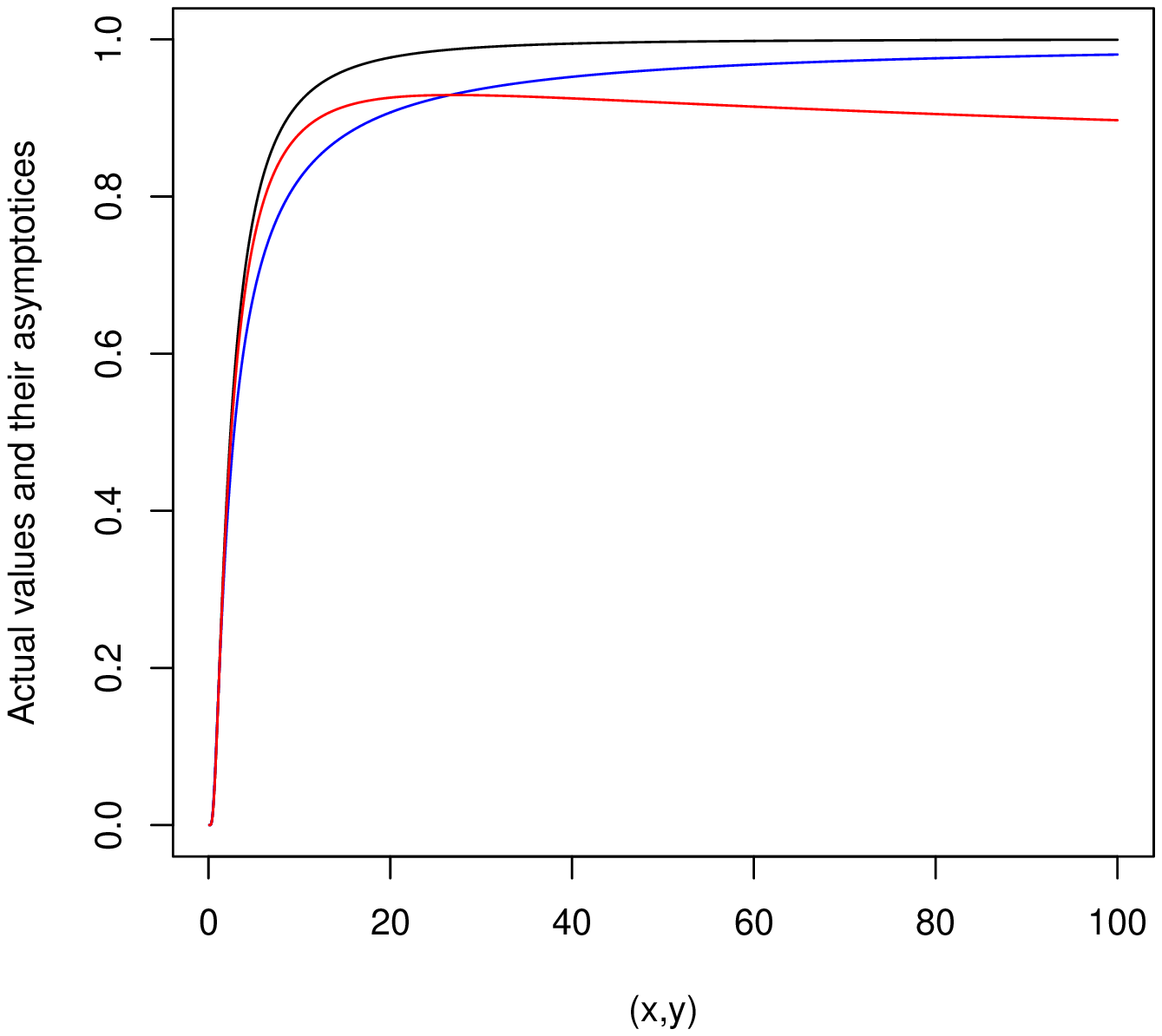, height=120pt, width=200pt,angle=0}
                   }%
                   \\
 \subfigure[$ \lambda=1,\tau=2, \rho_n=0.869$]
{%
\epsfig{file=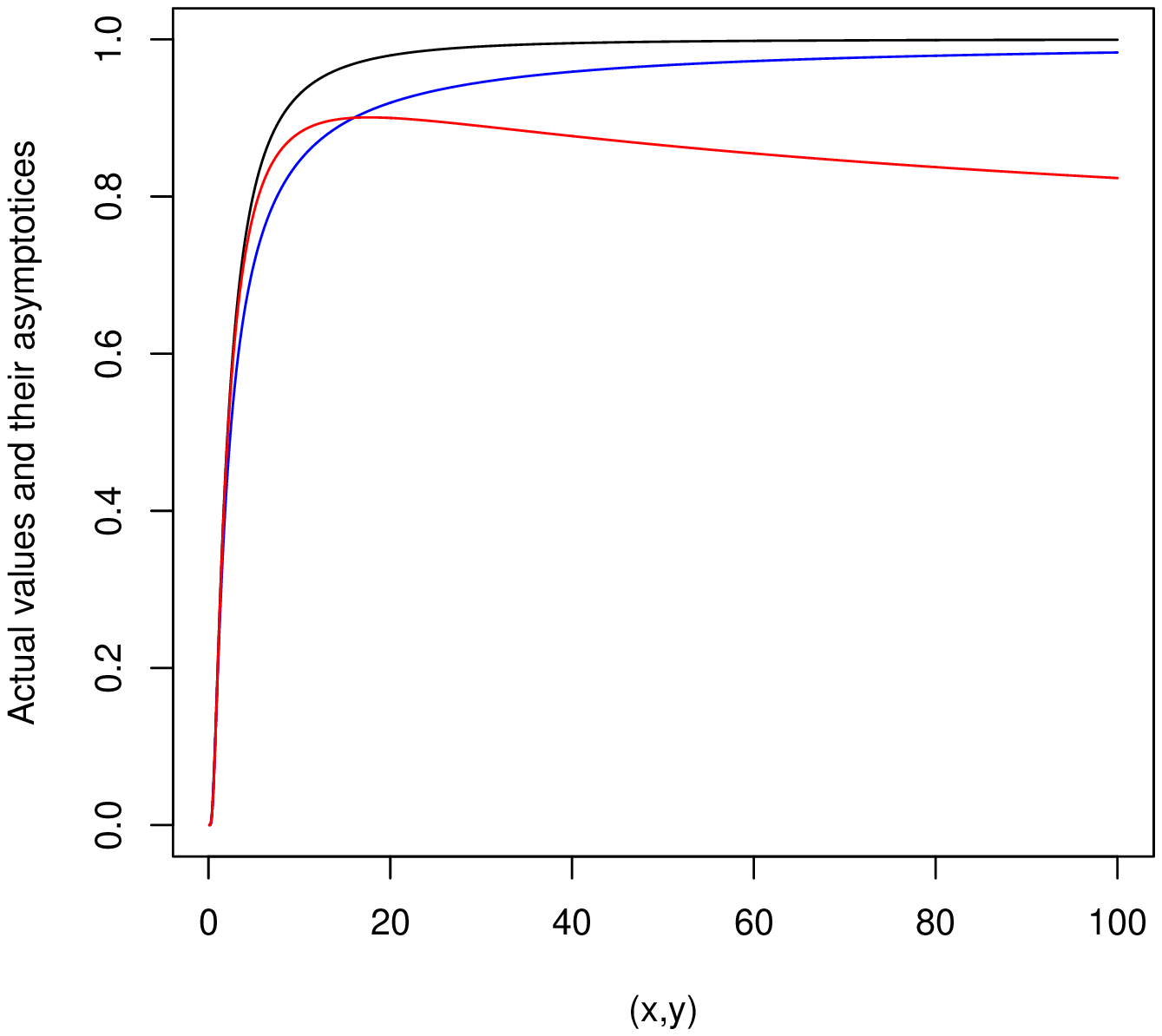, height=120pt, width=200pt,angle=0}
                   }%
\subfigure[$ \rho_n=1$]
{%
\epsfig{file=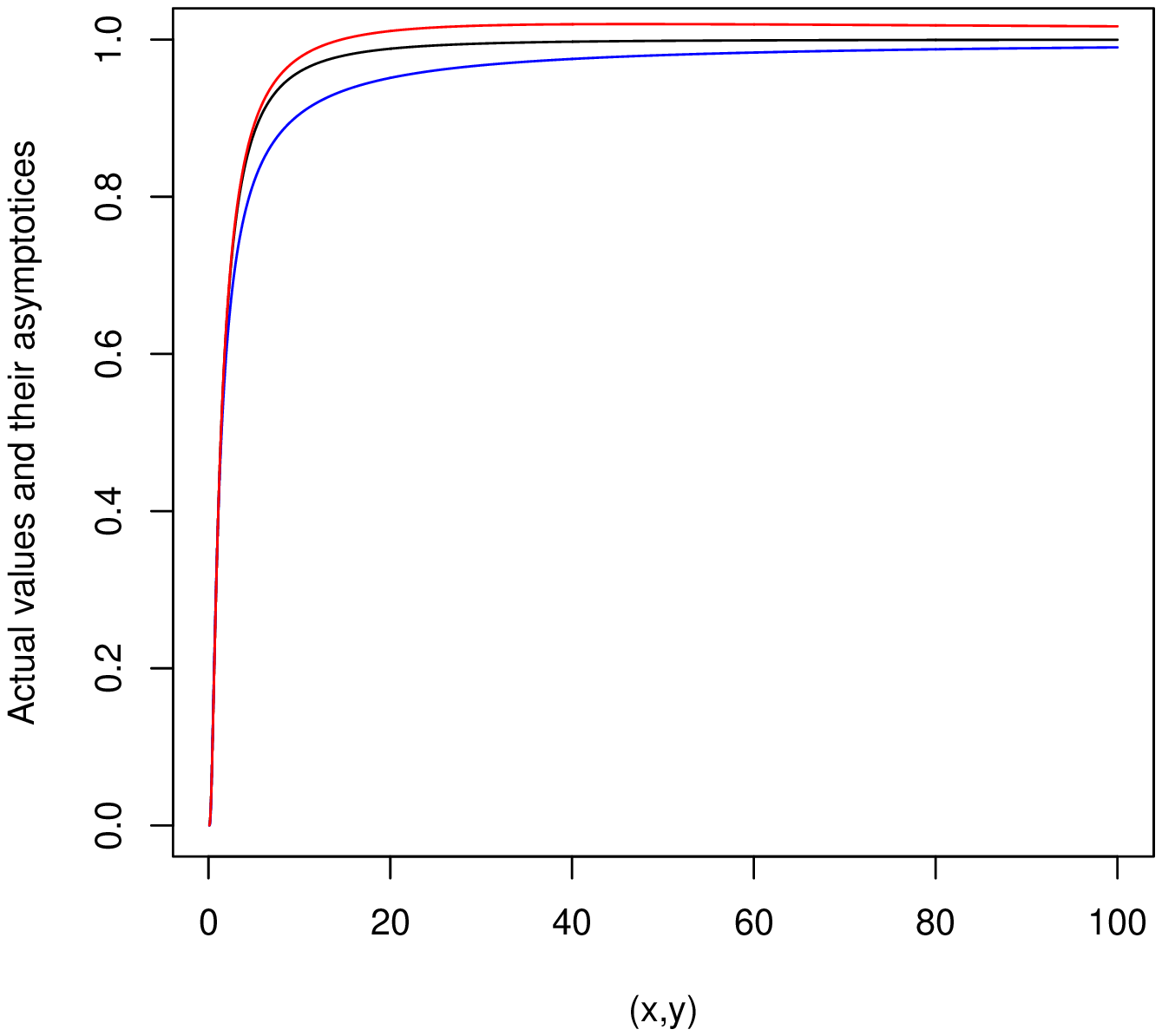, height=120pt, width=200pt,angle=0}
                   }%
\caption{Actual values and its approximations with $n=10^{3},
x=y\in [0,100]$. The actual values with black color, the first-order
asymptotics with blue color, the second-order asymptotics with red
color.}
\end{center}
\end{figure}

\begin{figure}\label{fig2}
\begin{center}
\subfigure[$\rho_n= -1$]
{%
 \epsfig{file=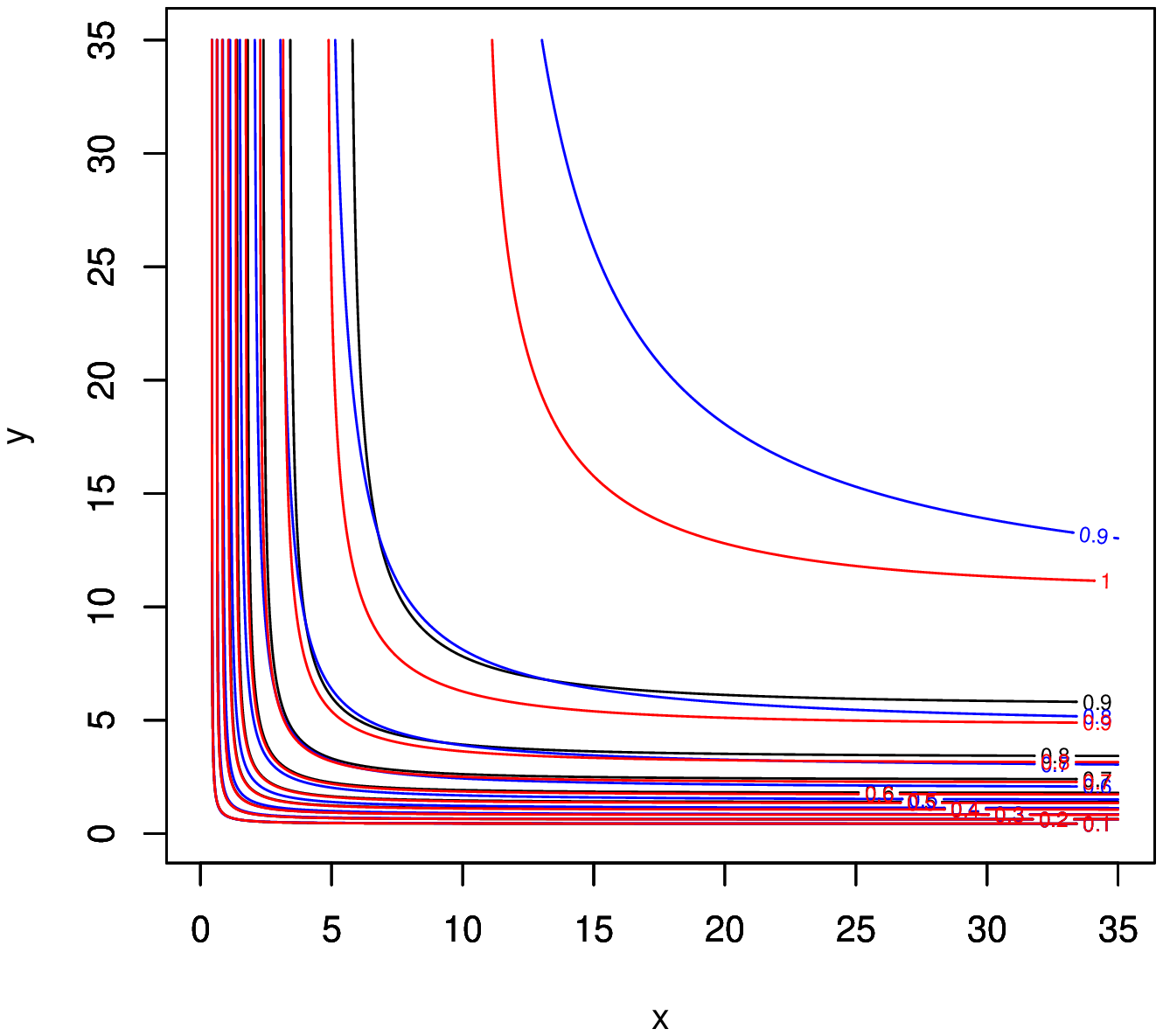, height=120pt, width=200pt,angle=0}
                   }%
\subfigure[$\lambda=2.5,\tau=-5, \rho_n=-0.915$]
{%
 \epsfig{file=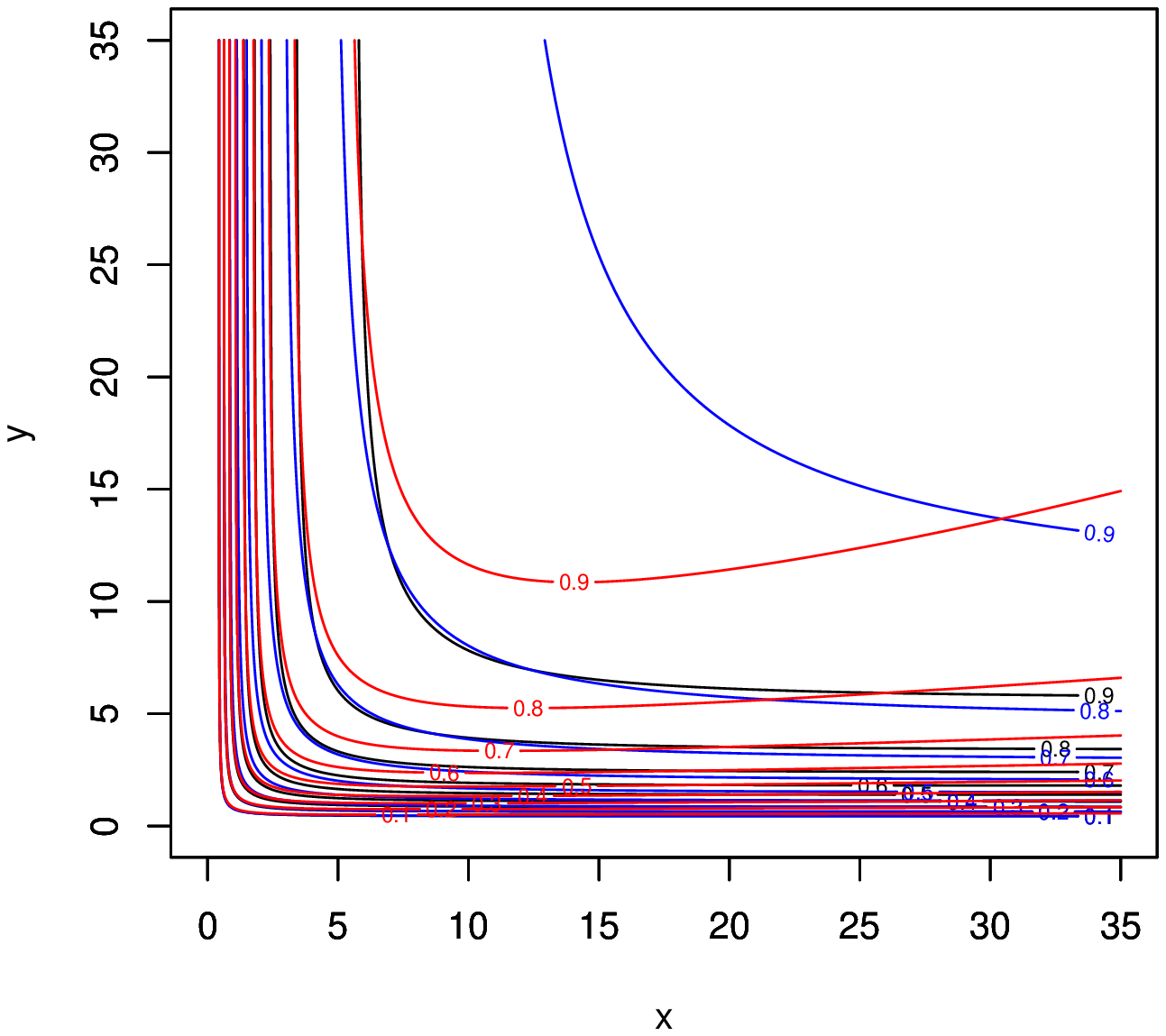, height=120pt, width=200pt,angle=0}
                   }%
                   \\
\subfigure[$\lambda=2.5,\tau=-2, \rho_n=-0.537$]
{%
 \epsfig{file=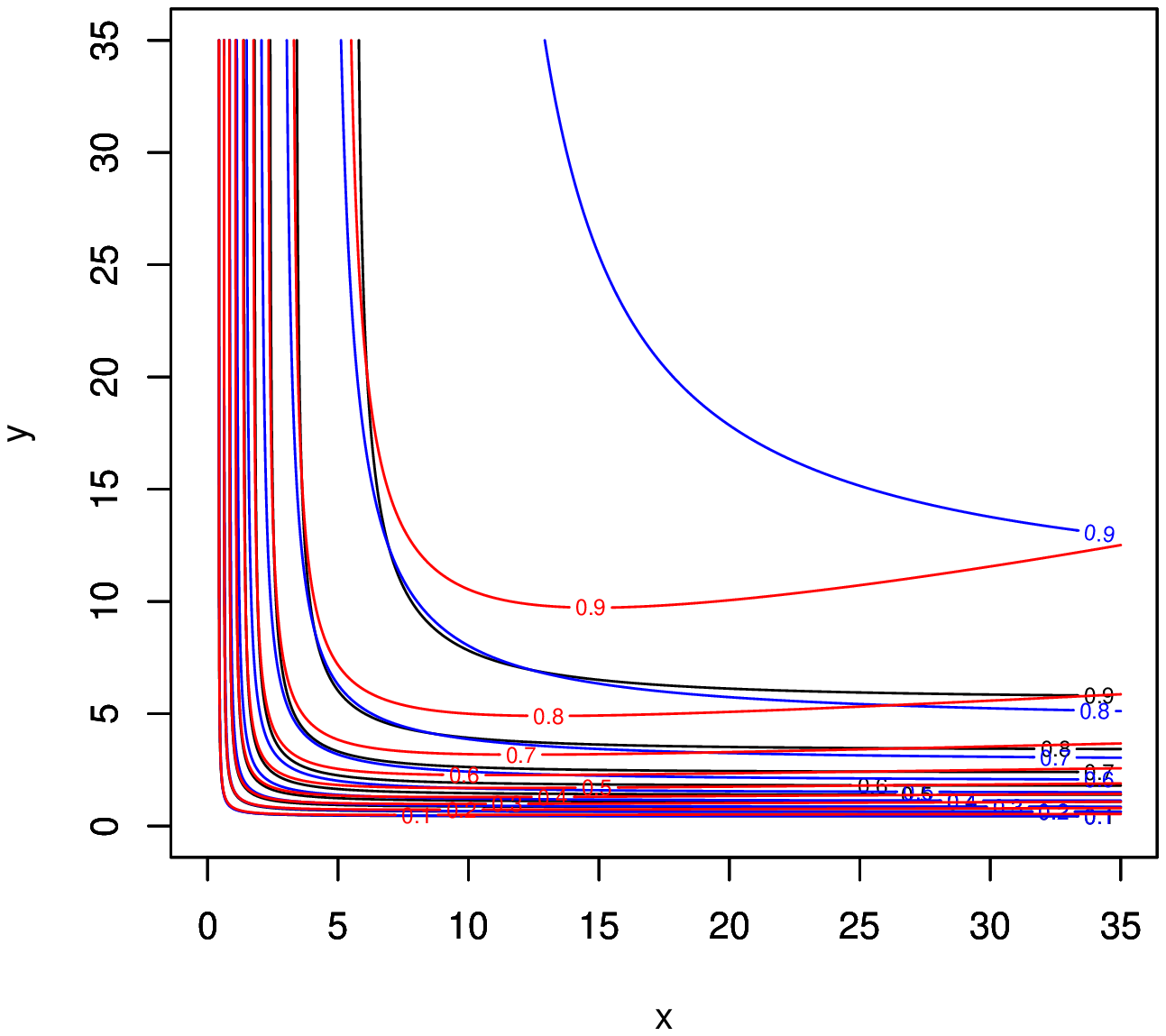, height=120pt, width=200pt,angle=0}
                   }%
 \subfigure[$\rho_n= 0$]
{%
 \epsfig{file=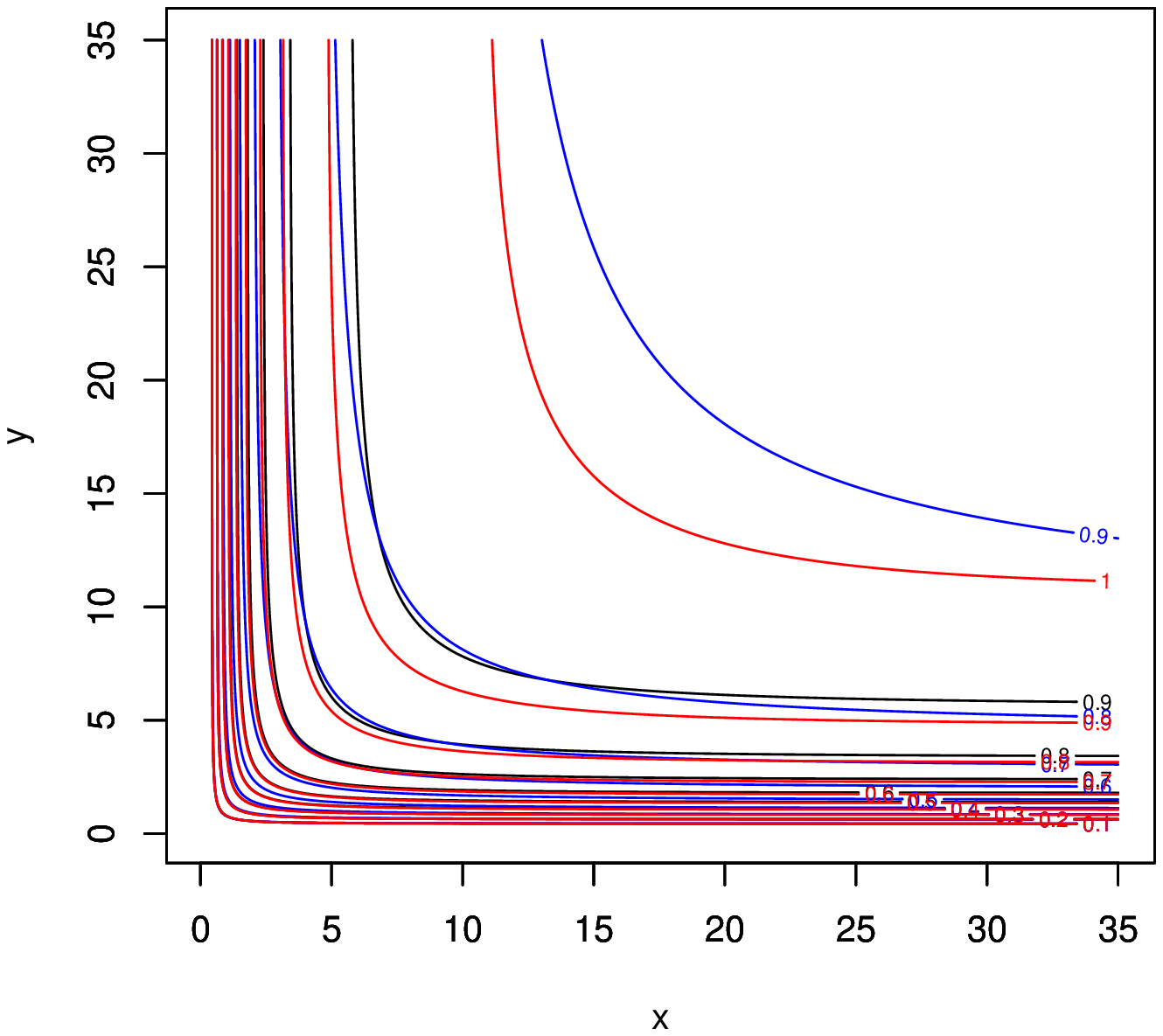, height=120pt, width=200pt,angle=0}
                   }%
                   \\
 \subfigure[$\lambda=2,\tau=2, \rho_n=0.329$]
{%
 \epsfig{file=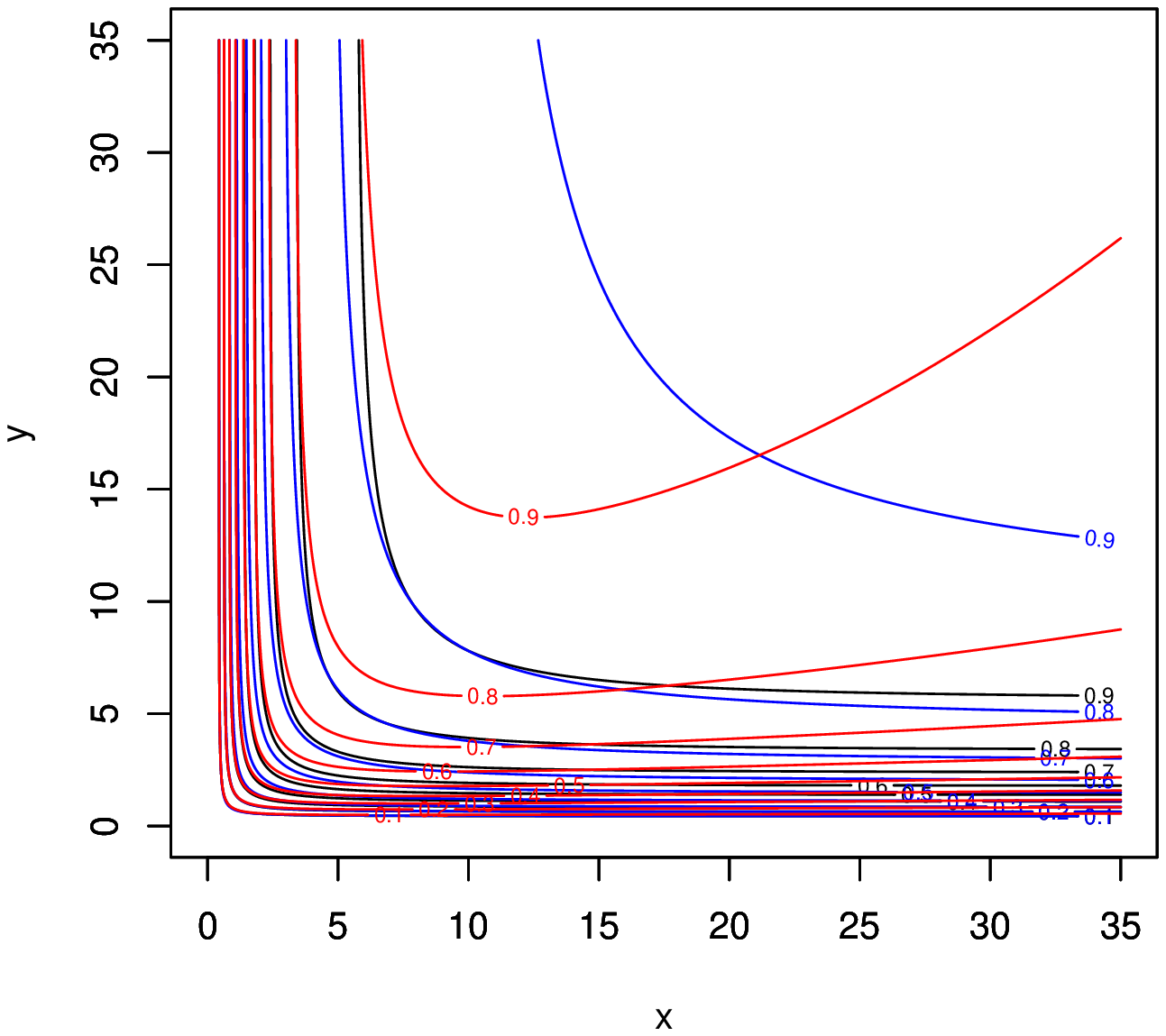, height=120pt, width=200pt,angle=0}
                   }%
\subfigure[$ \lambda=2,\tau=3,\rho_n=0.405$]
{%
\epsfig{file=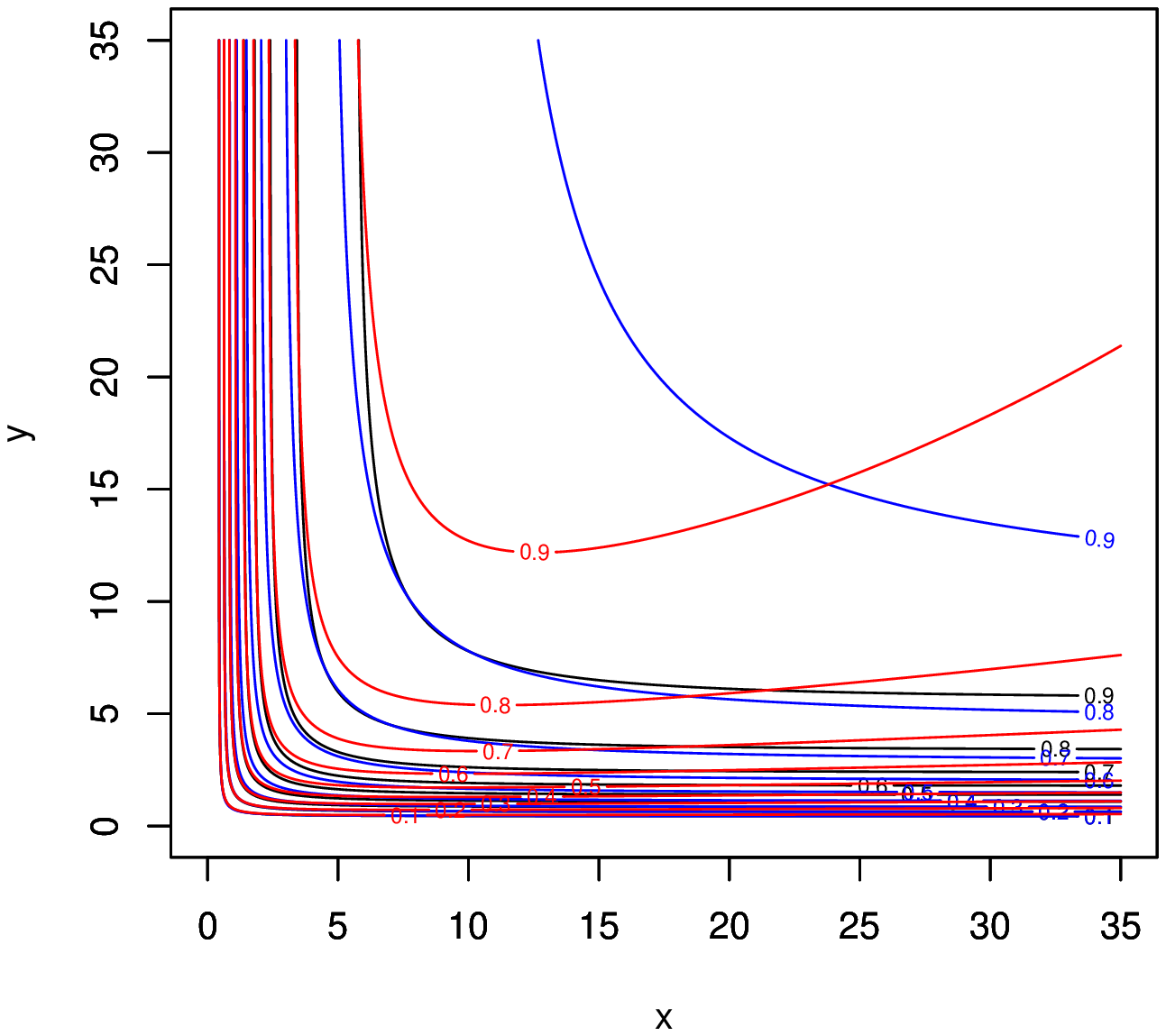, height=120pt, width=200pt,angle=0}
                   }%
                   \\
 \subfigure[$ \lambda=1,\tau=2, \rho_n=0.869$]
{%
\epsfig{file=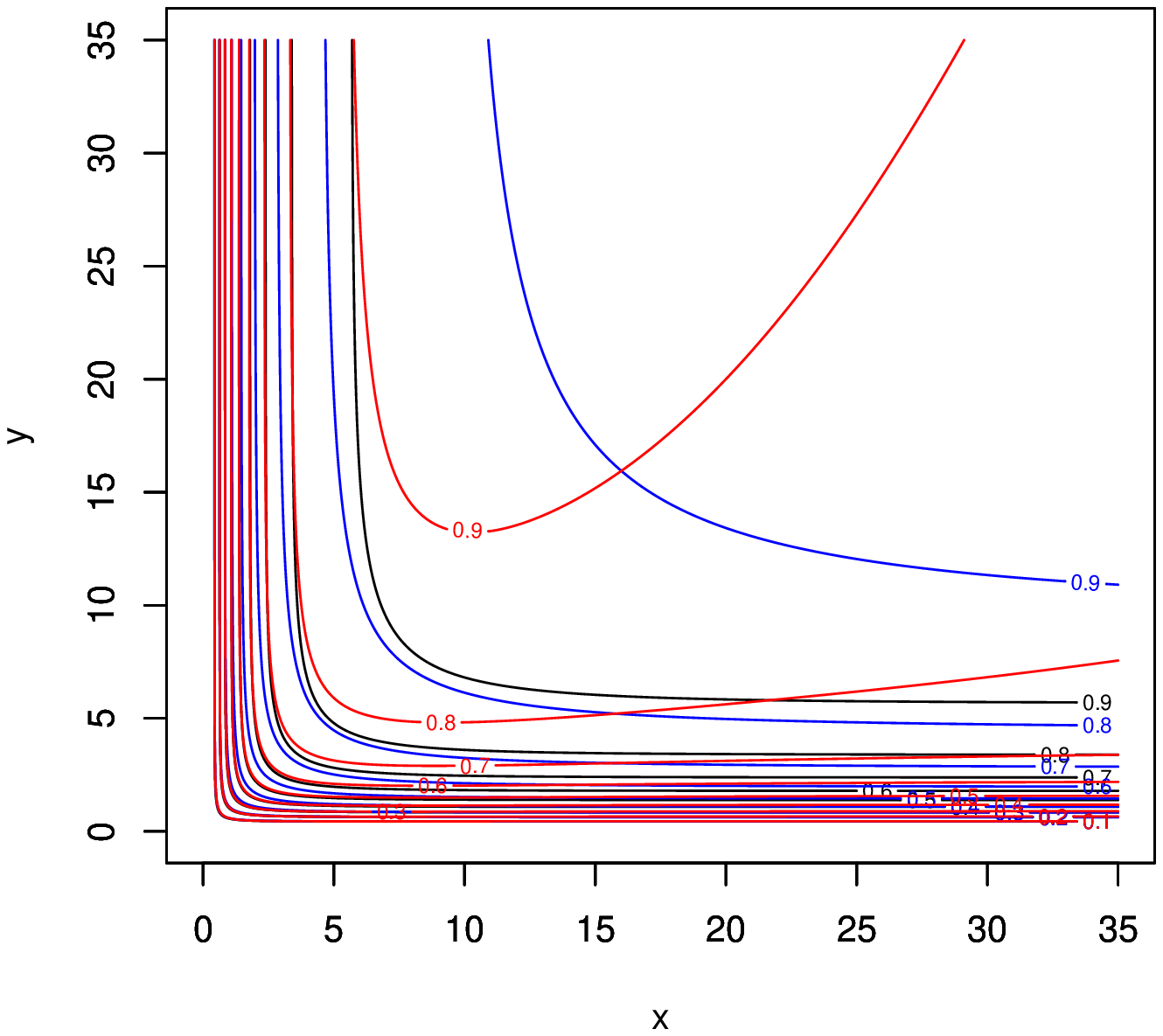, height=120pt, width=200pt,angle=0}
                   }%
\subfigure[$ \rho_n=1$]
{%
\epsfig{file=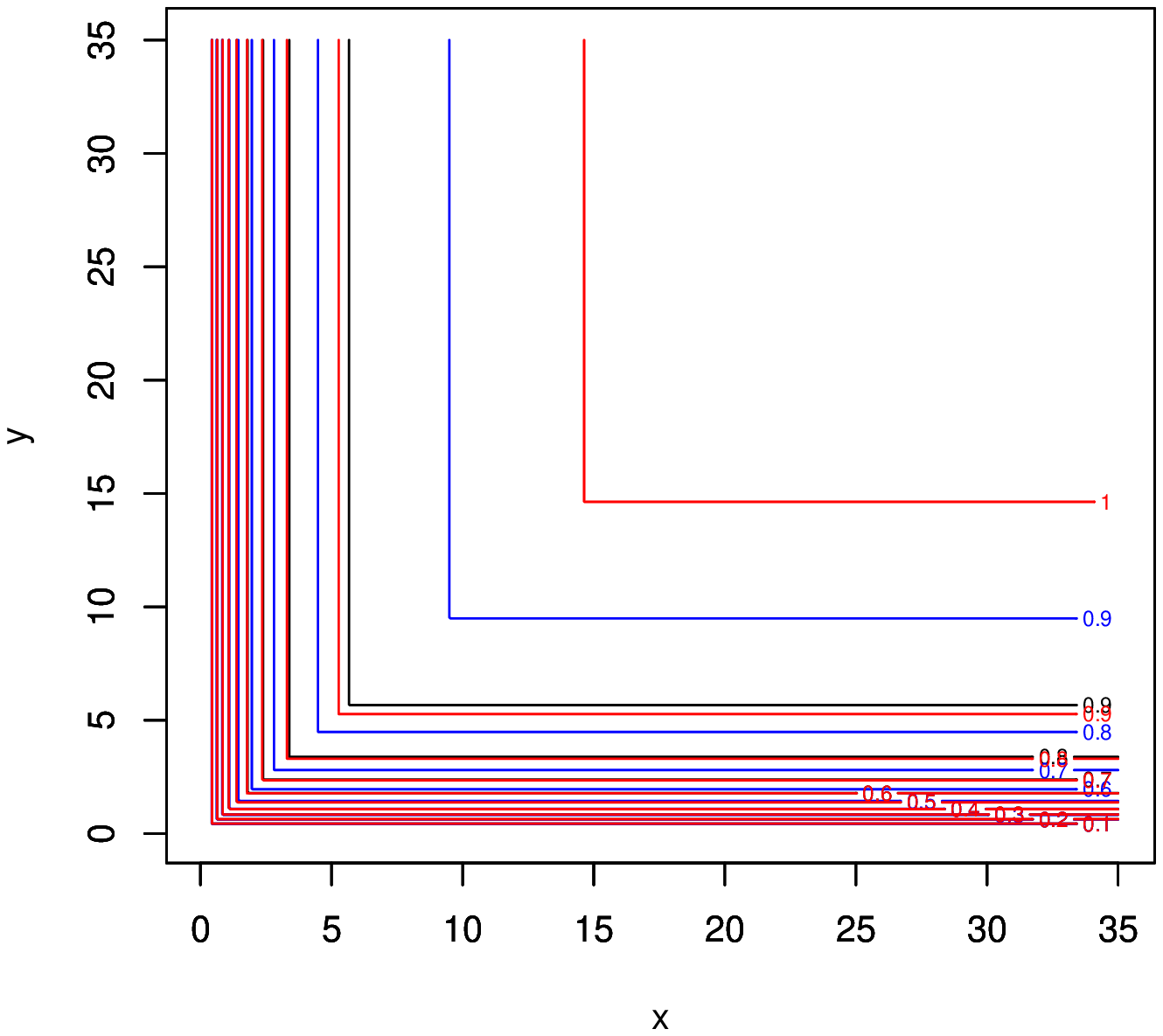, height=120pt, width=200pt,angle=0}
                   }%
\caption{The contour line of actual values and its approximations with $n=10^{3},
x,y\in [0,35]$. The actual values with black color, the first-order
asymptotics with blue color, the second-order asymptotics with red
color.}
\end{center}
\end{figure}

\section{Proofs}\label{sec3}

The aim of this section is to prove our main results.
Hereafter, for notational simplicity we shall write $u_n(x)=b_nx^{b_n^{-2}}$, $x>0$
with norming constant $b_n$ satisfying equation \eqref{eq1.2}.

\prooftheo{thm2.0}~The proofs of the theorem are similar with Lemma 21 in Kabluchko et al. (2009), so we omit here.
\qed

In order to prove Theorem \ref{th2.1}-Theorem \ref{th2.3}, we need some auxiliary lemmas as follows.
The following lemma shows the second order distributional expansions of maxima of
univariate Gaussian random sequences under power normalization, which proof is similar with Theorem 2.1 of Nair (1981).

\begin{lemma}\label{le3.1}
Let norming constants $b_n$ be satisfied \eqref{eq1.2}, then for $x>0$
$$\lim_{n\to \infty}b_n^2\Big(\Phi^n(u_n(x))-\exp\left(-x^{-1}\right)\Big)=s(x)\exp\left(-x^{-1}\right)$$
with
\begin{equation}\label{eq3.1}
s(x)=((\ln x)^2+\ln x)x^{-1}.
\end{equation}
\end{lemma}

\prooflem{le3.1}~ According to the definition of $b_n$ we have
\begin{equation}\label{eq3.2}
n^{-1}=1-\Phi(b_n)=b_n^{-1}\varphi(b_n)(1-b_n^{-2}+O(b_n^{-4}))
\end{equation}
with $\varphi(x)=\Phi^{\prime}(x)$ for large $n$, cf. Canto e Castro (1987). For $x>0$,
$n$ large, note
\begin{eqnarray*}
\frac{\varphi(u_n(x))}{\varphi(b_n)}&=&\exp\left(-\frac{1}{2}b_n^2(x^{2b_n^{-2}}-1)\right)\\
&=&\exp\left(-\frac{1}{2}b_n^2(e^{2b_n^{-2}\ln x}-1)\right)\\
&=&\exp\left(-\frac{1}{2}b_n^2(2b_n^{-2}\ln x+2b_n^{-4}(\ln x)^2+O(b_n^{-6}))\right)\\
&=&x^{-1}\exp\left(-b_n^{-2}(\ln x)^2(1+O(b_n^{-2}))\right)\\
&=&x^{-1}(1-(\ln x)^2b_n^{-2}+O(b_n^{-4})),
\end{eqnarray*}
we have
\begin{eqnarray*}
1-\Phi(u_n(x))
&=&b_n^{-1}x^{-b_n^{-2}}\varphi(b_n)\frac{\varphi(u_n(x))}{\varphi(b_n)}(1-b_n^{-2}x^{-2b_n^{-2}}+O(b_n^{-4}))\\
&=&b_n^{-1}\varphi(b_n)e^{-b_n^{-2}\ln x}x^{-1}\Big(1-(\ln x)^2b_n^{-2}+O(b_n^{-4})\Big)\Big(1-b_n^{-2}e^{-2b_n^{-2}\ln x}+O(b_n^{-4})\Big)\\
&=&b_n^{-1}\varphi(b_n)x^{-1}\Big(1-b_n^{-2}\ln x+O(b_n^{-4})\Big)\Big(1-(\ln x)^2b_n^{-2}+O(b_n^{-4})\Big)\Big(1-b_n^{-2}+O(b_n^{-4})\Big)\\
&=&b_n^{-1}\varphi(b_n)x^{-1}\Big(1-(1+\ln x+(\ln x)^2)b_n^{-2}+O(b_n^{-4})\Big).
\end{eqnarray*}
Let $h_n(x)=n \ln \Phi(u_n(x))+x^{-1}$, then
\begin{eqnarray*}
b_n^2h_n(x)&=&b_n^2\left(n\ln \Phi(u_n(x))+x^{-1}\right)\\
&=&b_n^2\Big(-n(1-\Phi(u_n(x)))-\frac{n}{2}(1-\Phi(u_n(x)))^2(1+o(1))+x^{-1}\Big)\\
&=&b_n^2\left[-x^{-1}\Big(1-(1+\ln x+(\ln x)^2)b_n^{-2}+O(b_n^{-4})\Big)(1+b_n^{-2}+O(b_n^{-4}))
-\frac{n}{2}(1-\Phi(u_n(x)))^2(1+o(1))+x^{-1}\right]\\
&\to&x^{-1}(\ln x+(\ln x)^2)=s(x),  \quad \mbox{as} \quad n\to \infty.
\end{eqnarray*}
Obviously, $\lim_{n\to \infty}h_n(x)=0$, thus
\begin{eqnarray*}
b_n^2\left(\Phi^n(u_n(x))-\exp\left(-x^{-1}\right)\right)
&=&b_n^2\left(\exp\left(h_n(x)\right)-1\right)\exp\left(-x^{-1}\right)\\
&=&b_n^2h_n(x)(1+o(1))\exp\left(-x^{-1}\right)\\
&\to&s(x)\exp\left(-x^{-1}\right)
\end{eqnarray*}
as $n\to \infty$. The proof is complete. \qed

The following two lemmas are mainly used to prove Theorem \ref{th2.1}.
A decomposition of probability $P(X>u_n(x),Y>u_n(y))$
is derived by Lemma \ref{le3.2}.
Lemma \ref{le3.3} gives the second order expansion of integration
$\int_{y}^{b_n^4}\Phi\left( \frac{u_n(x)-\rho_n u_n(z)}{\sqrt{1-\rho_n^2}} \right)z^{-2}dz$ by using refined condition \eqref{eq2.1} and Taylor expansion.
\begin{lemma}\label{le3.2}
If $(X,Y)$ be a bivariate Gaussian vector with correlation $\rho_{n} \in (-1,1)$, then
for $x,y>0$,
\begin{eqnarray}\label{addeq3.3}
&&n\pk{X>u_n(x),Y>u_n(y)} \nonumber \\
&=&n(1-\Phi(u_n(y)))
-\int_y^{b_{n}^{4}}\Phi\left(\frac{u_n(x)-\rho_{n} u_n(z)}{\sqrt{1-\rho_{n}^2}}\right)z^{-2}
\left[1+(1+\ln z-(\ln z)^2)b_n^{-2}\right]dz+O(b_n^{-4})
\end{eqnarray}
for large $n$.
\end{lemma}

\prooflem{le3.2}~First note that
\begin{eqnarray}\label{waddeq1}
1+x+\frac{x^{2}}{2}+\frac{x^{3}}{6}<e^{x}<1+x+\frac{x^{2}}{2}+\frac{x^{3}}{6}+x^4, 0<x<1
\end{eqnarray}
and
\begin{eqnarray}\label{waddeq2}
1-x+\frac{x^2}{2}-\frac{x^3}{6}<e^{-x}<1-x+\frac{x^2}{2}+\frac{x^3}{6}, x>0.
\end{eqnarray}
So we have
\begin{eqnarray}\label{addeq3.4}
& & \int_{u_{n}(y)}^{u_{n}(b_{n}^{4})}\Phi\left( \frac{u_{n}(x)-\rho_{n}z}{\sqrt{1-\rho_{n}^{2}}}¡¡\right)\varphi(z)dz \nonumber \\
&=& \int_{y}^{b_{n}^{4}} \Phi\left(\frac{u_n(x)-\rho_{n} u_n(z)}{\sqrt{1-\rho_{n}^2}}\right)
\varphi(u_n(z))b_n^{-1}z^{b_n^{-2}-1}dz \nonumber \\
&=&b_n^{-1}\varphi(b_n)\int_y^{b_{n}^{4}}\Phi\left(\frac{u_n(x)-\rho_{n} u_n(z)}{\sqrt{1-\rho_{n}^2}}\right)
e^{-\frac{1}{2}b_n^2(z^{2b_n^{-2}}-1)}e^{b_n^{-2}\ln z}z^{-1} dz \nonumber\\
&=&b_n^{-1}\varphi(b_n)\int_y^{b_{n}^{4}}\Phi\left(\frac{u_n(x)-\rho_{n} u_n(z)}{\sqrt{1-\rho_{n}^2}}\right)
\Big( 1- b_{n}^{-2}(\ln z)^{2} \Big)
\Big( 1+b_{n}^{-2}\ln z  \Big) z^{-2}dz+O(b_n^{-5}\varphi(b_n)) \nonumber \\
&=& b_n^{-1}\varphi(b_n)\int_y^{b_{n}^{4}}\Phi\left(\frac{u_n(x)-\rho_{n} u_n(z)}{\sqrt{1-\rho_{n}^2}}\right)
\Big( 1+ (\ln z -(\ln z)^{2} )b_{n}^{-2}  \Big) z^{-2} dz +O(b_n^{-5}\varphi(b_n))
\end{eqnarray}
for large $n$.

It is well known that
\begin{equation}\label{addeq3.5}
1-\Phi(x)<x^{-1}\varphi(x)
\end{equation}
for $x>0$.  Combining with the inequality $e^{x}\geq 1+x, x\in \R$,
we can get
\begin{eqnarray}\label{addeq3.6}
 \int_{u_{n}(b_{n}^{4})}^{\infty} \Phi\left( \frac{u_{n}(x)-\rho_{n}z}{\sqrt{1-\rho_{n}^{2}}}¡¡\right)\varphi(z)dz
&\leq& 1-\Phi\left( u_{n}\left(b_{n}^{4}\right) \right) \nonumber\\
&\leq& b_{n}^{-4b_{n}^{-2}-1}\varphi\left( b_{n}^{4b_{n}^{-2}+1} \right) \nonumber\\
&=& b_{n}^{-4b_{n}^{-2}-1}\varphi(b_{n})\exp\left( -\frac{b_{n}^{2}}{2}\left( b_{n}^{8b_{n}^{-2}} -1 \right) \right) \nonumber \\
&<& b_{n}^{-4b_{n}^{-2}-5}\varphi(b_{n}) \nonumber\\
&=&O\left( b_{n}^{-5}\varphi(b_{n}) \right).
\end{eqnarray}

Since
\begin{equation*}
n\pk{X>u_n(x),Y>u_n(y)}=n(1-\Phi(u_n(y)))-n\int_{u_n(y)}^{\infty}\Phi\left(\frac{u_n(x)-\rho_{n} z}{
\sqrt{1-\rho_{n}^2}}\right)\varphi(z)dz,
\end{equation*}
we can derive \eqref{addeq3.3} by combining with \eqref{eq3.2}, \eqref{addeq3.4}, \eqref{addeq3.6}.

The proof is complete.
 \qed

\begin{lemma}\label{le3.3}
Under the conditions of Theorem \ref{th2.1}, we have
$$\lim_{n\to\infty}b_n^2\int_y^{b_{n}^{4}}
\left(\Phi\left(\lambda+\frac{\ln \frac{x}{z}}{2\lambda}\right)-\Phi\left(\frac{u_n(x)-\rho_n u_n(z)}{\sqrt{1-\rho_n^2}}\right)\right)z^{-2}dz
=\kappa_1(x,y,\lambda,\tau),$$
where
\begin{eqnarray*}
\kappa_1(x,y,\lambda,\tau)&=&(4\lambda^4+2\lambda^2-4\lambda^2\ln x)x^{-1}\left(1-\Phi\left(\lambda+\frac{\ln \frac{y}{x}}{2\lambda}\right)\right)+(2\tau -5\lambda^3+\lambda \ln y+\lambda \ln x) x^{-1}\varphi\left(\lambda+\frac{\ln \frac{y}{x}}{2\lambda}\right).
\end{eqnarray*}
\end{lemma}

\prooflem{le3.3}~Using the assumption \eqref{eq2.1} we can get
$$\lim_{n\to \infty}b_n^2\left(\lambda-\lambda_n\left(1-\frac{\lambda_n^2}{b_n^2}\right)^{-\frac{1}{2}} \right)
=\tau-\frac{1}{2}\lambda^3,$$
$$\lim_{n\to \infty}\frac{1}{2}b_n^2\left(\frac{1}{\lambda}-\frac{1}{\lambda_n}
\left(1-\frac{\lambda_n^2}{b_n^2}\right)^{-\frac{1}{2}}\right)=-\frac{1}{2}\tau
\lambda^{-2}-\frac{1}{4}\lambda$$
and
$$\lim_{n\to \infty}\lambda_n\left(1-\frac{\lambda_n^2}{b_n^2}\right)^{-\frac{1}{2}}=\lambda.$$
Further, by partial integration we get
$$\int_y^{\infty}\varphi\left(\lambda+\frac{\ln \frac{x}{z}}{2\lambda}\right)z^{-2}dz
=2\lambda x^{-1}\left(1-\Phi\left(\lambda+\frac{\ln \frac{y}{x}}{2\lambda}\right)\right),$$
$$\int_y^{\infty}\varphi\left(\lambda+\frac{\ln \frac{x}{z}}{2\lambda}\right)z^{-2}\ln zdz
=4\lambda^2 x^{-1}\varphi\left(\lambda+\frac{\ln \frac{y}{x}}{2\lambda}\right)
+(2\lambda \ln x-4\lambda^3) x^{-1}\left(1-\Phi\left(\lambda+\frac{\ln \frac{y}{x}}{2\lambda}\right)\right)$$
and
\begin{eqnarray*}&&\int_y^{\infty}\varphi\left(\lambda+\frac{\ln \frac{x}{z}}{2\lambda}\right)z^{-2}(\ln z)^2dz\\
&=&(4\lambda^2x^{-1}\ln y+4\lambda^2 x^{-1}(\ln x-2\lambda^2))\varphi\left(\lambda+\frac{\ln \frac{y}{x}}{2\lambda}\right)
+(8\lambda^3+2\lambda(\ln x-2\lambda^2)^2) x^{-1}\left(1-\Phi\left(\lambda+\frac{\ln \frac{y}{x}}{2\lambda}\right)\right).
\end{eqnarray*}
\COM{
For $y<z<b_{n}^4$, by using \eqref{addeq3.13} we have
\begin{eqnarray*}
&&\left(\lambda_n+\frac{\ln x/z}{2\lambda_n}+\frac{\lambda_n\ln z}{b^2_n}+\frac{(\ln x)^2-(\ln z)^2}{4b_n^2\lambda_n}+\frac{(\ln z)^2\lambda_n}{2b_n^4}
+\frac{(\ln x)^3-\rho_n(\ln z)^3}{12b_n^4\lambda_n}-\frac{\rho_n (\ln z)^4}{2b_n^6\lambda_n}\right)\left(1-\frac{\lambda_n^2}{b_n^2}\right)^{-\frac{1}{2}}\\
&\le&\frac{u_n(x)-\rho_n u_n(z)}{\sqrt{1-\rho_n^2}}=\frac{b_ne^{b_n^{-2}\ln x}-\rho_nb_ne^{b_n^{-2}\ln z}}{\sqrt{1-\rho_n^2}}\\
&\le&\left(\lambda_n+\frac{\ln x/z}{2\lambda_n}+\frac{\lambda_n\ln z}{b^2_n}+\frac{(\ln x)^2-(\ln z)^2}{4b_n^2\lambda_n}+\frac{(\ln z)^2\lambda_n}{2b_n^4}
+\frac{(\ln x)^3-\rho_n(\ln z)^3}{12b_n^4\lambda_n} + \frac{(\ln x)^4}{2b_{n}^6\lambda_{n}}\right)\left(1-\frac{\lambda_n^2}{b_n^2}\right)^{-\frac{1}{2}}.
\end{eqnarray*}
}
Since
\begin{equation}\label{addeq3.8}
\frac{u_n(x)-\rho_n u_n(z)}{\sqrt{1-\rho_n^2}}\to \lambda+\frac{\ln \frac{x}{z}}{2\lambda},\quad n\to \infty,
\end{equation}
and using \eqref{waddeq1} and \eqref{waddeq2}, we have
\begin{eqnarray}\label{addeq3.9}
&&b_n^2\int_y^{b_{n}^{4}}\left(\lambda+\frac{\ln \frac{x}{z}}{2\lambda}-\frac{u_n(x)-\rho_n u_n(z)}{\sqrt{1-\rho_n^2}}\right)
\varphi\left(\lambda+\frac{\ln \frac{x}{z}}{2\lambda}\right)z^{-2}dz \nonumber\\
&=&b_n^2\int_y^{b_{n}^{4}}\left[\left(\lambda-\lambda_n \left(1-\frac{\lambda_n^2}{b_n^2}\right)^{-\frac{1}{2}}\right)
+\left(\ln \frac{x}{z}\right)\left(\frac{1}{2\lambda}-\frac{1}{2\lambda_n}\left(1-\frac{\lambda_n^2}{b_n^2}\right)^{-\frac{1}{2}}\right)\right.\nonumber\\
&&\qquad  \left.-\frac{\lambda_n\ln z}{b_n^2}\left(1-\frac{\lambda_n^2}{b_n^2}\right)^{-\frac{1}{2}}
-\frac{(\ln x)^2-(\ln z)^2}{4b_n^2\lambda_n}\left(1-\frac{\lambda_n^2}{b_n^2}\right)^{-\frac{1}{2}}\right]
\varphi\left(\lambda+\frac{\ln \frac{x}{z}}{2\lambda}\right)z^{-2}dz +O(b_n^{-2}) \nonumber\\
&\to&\left(\tau-\frac{1}{2}\lambda^3-\frac{1}{2}\tau \lambda^{-2}\ln x-\frac{1}{4}\lambda \ln x -\frac{1}{4\lambda}(\ln x)^2\right)
\int_y^{\infty}\varphi\left(\lambda+\frac{\ln \frac{x}{z}}{2\lambda}\right)z^{-2}dz\nonumber\\
&&+(\frac{1}{2}\tau \lambda^{-2}-\frac{3}{4}\lambda)
\int_y^{\infty}\varphi\left(\lambda+\frac{\ln \frac{x}{z}}{2\lambda}\right)z^{-2}\ln z dz
+\frac{1}{4\lambda}\int_y^{\infty}\varphi\left(\lambda+\frac{\ln \frac{x}{z}}{2\lambda}\right)z^{-2}(\ln z)^2 dz\nonumber\\
&=&(4\lambda^4+2\lambda^2-4\lambda^2\ln x)x^{-1}\left(1-\Phi\left(\lambda+\frac{\ln \frac{y}{x}}{2\lambda}\right)\right) +(2\tau -5\lambda^3+\lambda \ln y+\lambda \ln x )x^{-1}\varphi\left(\lambda+\frac{\ln \frac{y}{x}}{2\lambda}\right)\nonumber\\
&=:&\kappa_1(x,y,\lambda,\tau),
\end{eqnarray}
as $n\to \infty$.
Using Taylor expansion with Lagrange remainder term, we have
\begin{eqnarray}\label{addeq3.10}
\Phi\left(\frac{u_n(x)-\rho_n u_n(z)}{\sqrt{1-\rho_n^2}}\right)&=&
\Phi\left(\lambda+\frac{\ln \frac{x}{z}}{2\lambda}\right)
+\varphi\left(\lambda+\frac{\ln \frac{x}{z}}{2\lambda}\right)
\left(\frac{u_n(x)-\rho_n u_n(z)}{\sqrt{1-\rho_n^2}}-\lambda-\frac{\ln \frac{x}{z}}{2\lambda}\right)\nonumber\\
&&-\frac{1}{2}\xi_n\varphi(\xi_n)\left(\frac{u_n(x)-\rho_n u_n(z)}{\sqrt{1-\rho_n^2}}-\lambda-\frac{\ln \frac{x}{z}}{2\lambda}\right)^2,
\end{eqnarray}
for some $\xi_n$ between $\frac{u_n(x)-\rho_n u_n(z)}{\sqrt{1-\rho_n^2}}$ and $\lambda+\frac{\ln x/z}{2\lambda}$. Moreover, using
dominated convergence theorem and \eqref{addeq3.8} and \eqref{addeq3.9},
\begin{equation}\label{addeq3.11}
b_n^2\int_y^{b_n^4}\left(\frac{u_n(x)-\rho_n u_n(z)}{\sqrt{1-\rho_n^2}}-\lambda-\frac{\ln \frac{x}{z}}{2\lambda}\right)^2
\xi_n\varphi(\xi_n)z^{-2}dz=O(b_n^{-2}).
\end{equation}
Combining with \eqref{addeq3.9}-\eqref{addeq3.11}, we get the desired result. \qed

With the above three lemmas, now we can give the proof of Theorem \ref{th2.1}.

\prooftheo{th2.1}~ Define
$$h_n(x,y,\lambda)=n\ln F(u_n(x),u_n(y))+\Phi\left(\lambda+\frac{\ln \frac{x}{y}}{2\lambda}\right)y^{-1}
+\Phi\left(\lambda+\frac{\ln \frac{y}{x}}{2\lambda}\right)x^{-1}.$$
In view of Lemma \ref{le3.1}, Lemma \ref{le3.2} and Lemma \ref{le3.3}, we have
\begin{eqnarray*}
b_n^2h_n(x,y,\lambda)&=&b_n^2\left[-n(1-F(u_n(x),u_n(y)))
-\frac{n}{2}(1-F(u_n(x),u_n(y)))^2(1+o(1))\right.\\
&&\qquad \left.+\Phi\left(\lambda+\frac{\ln \frac{x}{y}}{2\lambda}\right)y^{-1}
+\Phi\left(\lambda+\frac{\ln \frac{y}{x}}{2\lambda}\right)x^{-1}\right]\\
&=&b_n^2\left[-n(1-\Phi(u_n(x)))+x^{-1}\right]
+b_n^2\int_y^{b_{n}^{4}}\left[\Phi\left(\lambda+\frac{\ln \frac{x}{z}}{2\lambda}\right)-
\Phi\left(\frac{u_n(x)-\rho_n u_n(z)}{\sqrt{1-\rho_n^2}}\right)\right]z^{-2}dz\\
&&-\int_y^{b_{n}^{4}}\Phi\left(\frac{u_n(x)-\rho_n u_n(z)}{\sqrt{1-\rho_n^2}}\right)z^{-2}(1+\ln z-(\ln z)^2)dz
+O(b_n^{-2})-\frac{nb_n^2}{2}(1-F(u_n(x),u_n(y)))^2(1+o(1))\\
&\to&s(x)+\kappa_1(x,y,\lambda,\tau)
-\int_y^{\infty}\Phi\left( \lambda +\frac{\ln\frac{x}{z}}{2\lambda} \right)z^{-2}(1+\ln z-(\ln z)^2)dz
\end{eqnarray*}
as $n\to \infty$. By partial integration we have
\begin{eqnarray*}
&& \int_y^{\infty}\Phi\left( \lambda +\frac{\ln\frac{x}{z}}{2\lambda} \right)z^{-2}(1+\ln z-(\ln z)^2)dz\\
&=&\Big((\ln x)^2 + \ln x -4\lambda^2(\ln x) +2\lambda^2+4\lambda^4 \Big)x^{-1}
\left(1-\Phi\left(\lambda+\frac{\ln \frac{y}{x}}{2\lambda}\right)\right)\\
&&+\Big(2\lambda+2\lambda (\ln y) +2\lambda (\ln x) -4\lambda^3\Big)x^{-1}\varphi\left(\lambda+\frac{\ln \frac{y}{x}}{2\lambda}\right)\\
&&-\Big((\ln y)^2 + \ln y\Big)y^{-1}\Phi\left(\lambda+\frac{\ln \frac{x}{y}}{2\lambda}\right).
\end{eqnarray*}

Obviously, $h_n(x,y,\lambda)\to 0,$
hence
\begin{eqnarray*}
&&b_n^2(F^n(u_n(x),u_n(y))-H_{\lambda}(\ln x,\ln y))\\
&=&b_n^2(\exp\left(h_n(x,y,\lambda)\right)-1)H_{\lambda}(\ln x,\ln y)\\
&=&b_n^2h_n(x,y,\lambda)(1+o(1))H_{\lambda}(\ln x,\ln y)\\
&\to&\kappa(x,y,\lambda,\tau)H_{\lambda}(\ln x,\ln y)
\end{eqnarray*}
as $n\to \infty$, where
\begin{eqnarray}\label{addeq3.12}
\kappa(x,y,\lambda,\tau)
&=&\Big((\ln x)^2 + \ln x\Big)x^{-1}\Phi\left(\lambda+\frac{\ln \frac{y}{x}}{2\lambda}\right)
+\Big((\ln y)^2 + \ln y\Big)y^{-1}\Phi\left(\lambda+\frac{\ln \frac{x}{y}}{2\lambda}\right)\nonumber\\
&&+(2\tau -\lambda^3 -2\lambda -\lambda \ln y -\lambda \ln x)\varphi\left(\lambda+\frac{\ln y/x}{2\lambda}\right).
\end{eqnarray}
The proof is complete. \qed

Next we prove the results of two extreme cases. For the case of $\lambda=\infty$,
\eqref{eq2.2} are derived by discussing $\rho_n=-1$ and $\rho_n\in (-1,1)$, respectively.

\prooftheo{th2.2}~ Let $h_n(x,y)=n\ln F(u_n(x),u_n(y))+x^{-1}+y^{-1}$.
First, we consider that the bivariate Gaussian is complete negative dependent $\rho_n= -1$ which implies $\lambda=\infty$.

According to Lemma \ref{le3.1}, for $\rho_n=-1$ we can get
\begin{eqnarray*}
b_n^2h_n(x,y)&=&b_n^2(-n(1-\Phi(u_n(x)))+x^{-1})+b_n^2(-n(1-\Phi(u_n(y)))+y^{-1})\\
&& -\frac{1}{2}b_n^2n(1-F(u_n(x),u_n(y)))^2(1+o(1))\\
&\to&s(x)+s(y), \quad \mbox{as}\quad n\to \infty.
\end{eqnarray*}

Consequently,
\begin{equation*}
\lim_{n\to \infty}b_n^2(F^n(u_n(x),u_n(y))-H_{\infty}(\ln x,\ln y))=(s(x)+s(y))H_{\infty}(\ln x,\ln y)
\end{equation*}
holds with $\rho_n= -1$.

Next we prove \eqref{eq2.2} holds with $\rho_n\in(-1,1)$.
For $\rho_n \in (-1,1)$ and fixed $x,z>0$, we have
$$\frac{u_n(x)-\rho_nu_n(z)}{\sqrt{1-\rho_n^2}}>0$$
for large $n$, due to $\lim_{n\to \infty}b_n^2(1-\rho_n)=\infty$.

Combining \eqref{addeq3.5} and the condition
$\lim_{n\to \infty}\frac{\ln b_n}{b_n^2(1-\rho_n)}=0$, we can get
\begin{eqnarray*}
&&b_n^3\left(1-\Phi\left(\frac{u_n(x)-\rho_nu_n(z)}{\sqrt{1-\rho_n^2}}\right)\right)\\
&\le& \frac{b_n^3\varphi\left(\frac{u_n(x)-\rho_nu_n(z)}{\sqrt{1-\rho_n^2}}\right)}
{b_n\sqrt{\frac{1-\rho_n}{1+\rho_n}}\left(1+\frac{\ln \frac{x}{z}}{b_n^2 (1-\rho_n)}+\frac{\ln z}{b_n^2}
+ O\left( \frac{(\ln b_{n})^2}{b_{n}^4(1-\rho_{n})} \right)\right)}\\
&\le&\frac{b_n^2\sqrt{\frac{1+\rho_n}{1-\rho_n}}\exp\left(-\frac{b_n^2(1-\rho_n)}{2(1+\rho_n)}\left(
1+ \frac{2(\ln x-\rho_n\ln z)}{b_{n}^2(1-\rho_n)}+O\left( \frac{(\ln b_{n})^2}{b_{n}^4(1-\rho_{n})} \right)  \right)\right) }
{\sqrt{2\pi}\left(1+\frac{\ln \frac{x}{z}}{b_n^{2}(1-\rho_n)}+\frac{\ln z}{b_n^2}+ O\left( \frac{(\ln b_{n})^2}{b_{n}^4(1-\rho_{n})} \right)\right)}\\
&\le&\frac{\exp\left(-\frac{b_n^2(1-\rho_n)}{2(1+\rho_n)}\left(1+\frac{2(\ln x-\ln z)}{b_n^2(1-\rho_n)}
+ \frac{2\ln z}{b_n^2}
-\frac{6(1+\rho_n)\ln b_n}{b_n^2(1-\rho_n)}+\frac{(1+\rho_n)\ln (b_n^2(1-\rho_n))}{b_n^2(1-\rho_n)}
+O\left( \frac{(\ln b_{n})^2}{b_{n}^4(1-\rho_{n})} \right) \right)\right)}
{\sqrt{\pi}\left(1+\frac{\ln \frac{x}{z}}{b_n^{2}(1-\rho_n)}+\frac{\ln z}{b_n^2}
+O\left( \frac{(\ln b_{n})^2}{b_{n}^4(1-\rho_{n})} \right)\right)}\\
&\le& \frac{\exp\left(-\frac{b_n^2(1-\rho_n)}{2(1+\rho_n)}\left(1+\frac{2(\ln x-4\ln b_{n})}{b_n^2(1-\rho_n)}
+ \frac{2\ln y}{b_n^2}
-\frac{6(1+\rho_n)\ln b_n}{b_n^2(1-\rho_n)}+\frac{(1+\rho_n)\ln (b_n^2(1-\rho_n))}{b_n^2(1-\rho_n)}
+O\left( \frac{(\ln b_{n})^2}{b_{n}^4(1-\rho_{n})} \right) \right)\right)}
{\sqrt{\pi}\left(1+\frac{\ln x - 4\ln b_{n}}{b_n^{2}(1-\rho_n)}+\frac{\ln y}{b_n^2}
+O\left( \frac{(\ln b_{n})^2}{b_{n}^4(1-\rho_{n})} \right)\right)} \\
&\to& 0
\end{eqnarray*}
as $n\to \infty$, if $y<z<b_{n}^4$. Hence
\begin{eqnarray*}
&&\pk{X>u_n(x),Y>u_n(y)}\\
&=&n^{-1}(1-b_n^{-2}+O(b_n^{-4}))^{-1}b_n^{-3}\left(\int_y^{b_{n}^{4}}b_n^{3}
\left(1-\Phi\left(\frac{u_n(x)-\rho_nu_n(z)}{\sqrt{1-\rho_n^2}}\right)\right)z^{-2}(1+(\ln z-(\ln z)^2)b_n^{-2})dz+O(b_n^{-1}) \right)\\
&=&O(n^{-1}b_n^{-3}).
\end{eqnarray*}
Using Lemma \ref{le3.1},
\begin{eqnarray*}
b_n^2h_n(x,y)&=&b_n^2(-n(1-\Phi(u_n(x)))+x^{-1})+b_n^2(-n(1-\Phi(u_n(y)))+y^{-1})\\
&&+b_n^2n\pk{X>u_n(x),Y>u_n(y)}-\frac{1}{2}b_n^2n(1-F(u_n(x),u_n(y)))^2(1+o(1))\\
&\to&s(x)+s(y), \quad \mbox{as}\quad n\to \infty.
\end{eqnarray*}
Thus, the claimed result \eqref{eq2.2} holds for $\rho_n\in (0,1)$, which complete
the proof.
 \qed

Similar with the proof of Theorem \ref{th2.2}, we prove the result for the extreme case of $\lambda=0$ by considering $\rho_n=1$ and $\rho_n\in (0,1)$.

\prooftheo{th2.3}~ For the complete positive dependent case $\rho_n\equiv1$,
\eqref{eq2.3} can be derived by combining Lemma \ref{le3.1}, so the rest is for the case of $\rho_n \in (0,1)$.

For $\rho_n \in (0,1)$ and fixed $x,y>0$, we have
$$\frac{u_n(\min(x,y))-\rho_nu_n(z)}{\sqrt{1-\rho_n^2}}<0,$$
if $\max(x,y)<z< b_n^{4}$.
Combining with $\lim_{n\to \infty} b_n^2(1-\rho_n)=0$, we can get
\begin{eqnarray}\label{eq3.13}
& & \Phi\left(\frac{u_n(\min(x,y))-\rho_nu_n(z)}{\sqrt{1-\rho_n^2}}\right)\nonumber\\
&\le&-\frac{\varphi\left(\frac{u_n(\min(x,y))-\rho_nu_n(z)}{\sqrt{1-\rho_n^2}}\right)}
{\frac{u_n(\min(x,y))-\rho_nu_n(z)}{\sqrt{1-\rho_n^2}}}\nonumber\\
&\le&\frac{\exp\left(-\frac{b_n^2(1-\rho_n)}{2(1+\rho_n)} -\frac{\ln \min(x,y)-\rho_n\ln z}{1+\rho_n} +O\left( \frac{(\ln b_{n})^2}{b_{n}^2} \right) \right)}
{\sqrt{2\pi}\left(\frac{\ln z-\ln \min(x,y)}{b_n\sqrt{1-\rho_n^2}}-\frac{b_n\sqrt{1-\rho_n}}{\sqrt{1+\rho_n}}
-\frac{\sqrt{1-\rho_n}\ln z}{b_n \sqrt{1+\rho_n}} + O\left( \frac{(\ln b_{n})^2}{b_{n}^3\sqrt{1-\rho_{n}^2}}  \right)\right)}\nonumber\\
&\le& \frac{b_{n}\sqrt{1-\rho_{n}^2}\exp\left(-\frac{b_n^2(1-\rho_n)}{2(1+\rho_n)} -\frac{\ln \min(x,y)-\rho_n\ln z}{1+\rho_n} +O\left( \frac{(\ln b_{n})^2}{b_{n}^2} \right) \right)}
{\sqrt{2\pi}\left(  \ln \max(x,y)-\ln \min(x,y) - 4(1-\rho_{n})\ln b_{n}
-b_{n}^2(1-\rho_{n}) + O\left( \frac{(\ln b_{n})^2}{b_{n}^2}  \right)\right)},
\end{eqnarray}
for large $n$ due to $\Phi(-x)=1-\Phi(x)$ and Mills' inequality \eqref{addeq3.5}.

From \eqref{eq3.13} and the inequality $e^{x}\geq 1+x, x\in R$, it follows that
\begin{eqnarray}\label{eq3.14}
& & \int_{\max(x,y)}^{b_{n}^{4}} \Phi\left( \frac{u_{n}(\min(x,y))-\rho_{n}u_{n}(z)}{\sqrt{1-\rho_{n}^{2}}} \right)
\exp\left( -\frac{b_{n}^{2}}{2} (z^{2b_n^{-2}}-1)\right)z^{b_{n}^{-2}}z^{-1}dz \nonumber \\
&\leq& \int_{\max(x,y)}^{b_{n}^{4}} \Phi\left( \frac{u_{n}(\min(x,y))-\rho_{n}u_{n}(z)}{\sqrt{1-\rho_{n}^{2}}} \right)
\exp\left( -\frac{b_{n}^{2}}{2}(2b_{n}^{-2}\ln z) \right) z^{b_{n}^{-2}}z^{-1}dz\nonumber\\
&\le& \frac{b_{n}\sqrt{1-\rho_{n}^2}\exp\left(-\frac{b_n^2(1-\rho_n)}{2(1+\rho_n)} -\frac{\ln \min(x,y)}{1+\rho_n} +O\left( \frac{(\ln b_{n})^2}{b_{n}^2} \right) \right) \int_{\max(x,y)}^{b_{n}^{4}} z^{b_{n}^{-2}}z^{-\frac{2+\rho_{n}}{1+\rho_{n}}}dz}
{\sqrt{2\pi}\left(  \ln \max(x,y)-\ln \min(x,y) - 4(1-\rho_{n})\ln b_{n}
-b_{n}^2(1-\rho_{n}) + O\left( \frac{(\ln b_{n})^2}{b_{n}^2}  \right)\right)}
\nonumber \\
&<& \frac{2b_{n}^{4b_{n}^{-2}} b_{n}\sqrt{1-\rho_{n}}\exp\left(-\frac{b_n^2(1-\rho_n)}{2(1+\rho_n)}
-\frac{\ln \min(x,y) + \ln \max(x,y)}{1+\rho_n} +O\left( \frac{(\ln b_{n})^2}{b_{n}^2} \right) \right)  }
{\sqrt{\pi}\left(  \ln \max(x,y)-\ln \min(x,y) - 4(1-\rho_{n})\ln b_{n}
-b_{n}^2(1-\rho_{n}) + O\left( \frac{(\ln b_{n})^2}{b_{n}^2}  \right)\right)} \nonumber \\
&=& o\left( b_{n}^{-2} \right)
\end{eqnarray}
for large $n$ by using $\lim_{n\to \infty} b_{n}^{6}(1-\rho_{n})=0$.

Form the proof of Lemma \ref{le3.1}, it follows that
\[1-\Phi(u_{n}(x))=n^{-1}\Big( x^{-1} -b_{n}^{-2}s(x)(1+o(1)) \Big).\]
Combining with \eqref{addeq3.6}, \eqref{eq3.14}, we have
\begin{eqnarray}
& & 1-F(u_{n}(x),u_{n}(y))\nonumber \\
&=& 1-\Phi\left( u_{n}(\min(x,y)) \right) + \int_{u_{n}(\max(x,y))}^{\infty}
\Phi\left( \frac{u_{n}(\min(x,y))-\rho_{n}z}{\sqrt{1-\rho_{n}^{2}}}\right)\varphi(z)dz  \nonumber \\
&=& n^{-1}\Big( (\min(x,y))^{-1} -b_{n}^{-2}s(\min(x,y))(1+o(1))   \nonumber \\
& & + \left( 1-b_{n}^{-2} + O(b_{n}^{-4}) \right)^{-1}\int_{\max(x,y)}^{b_{n}^{4}}
\Phi\left( \frac{u_{n}(\min(x,y))-\rho_{n}u_{n}(z) }{\sqrt{1-\rho_{n}^{2}}} \right)
\exp\left( -\frac{b_{n}^{2}}{2}\left( z^{2b_{n}^{-2}} -1 \right) \right)z^{b_{n}^{-2}}z^{-1}dz +O(b_{n}^{-4})  \Big) \nonumber \\
&=& n^{-1}\Big( (\min(x,y))^{-1} -b_{n}^{-2}s(\min(x,y))(1+o(1)) + o(b_{n}^{-2})  \Big)
\end{eqnarray}
for large $n$, which implies the desired result.

The proof is complete. \qed

\end{document}